\documentclass[11pt]{article}
\usepackage{amsmath,amsfonts,amsthm,amssymb,amscd}
\usepackage[latin1]{inputenc}
\usepackage{appendix}
\usepackage{latexsym}
\usepackage{graphicx}
\usepackage{enumerate}

\def\M{\mathbb{N}_{0}}
\def\N{\mathbb{N}}

\def\m{\text{M}}
\def\w{\text{W}}
\def\P{\mathcal{P}}

\newtheorem{Thm}{Theorem}[section]
\newtheorem{Prop}[Thm]{Proposition}
\newtheorem{Lem}[Thm]{Lemma}

\newtheorem{Conj}[Thm]{Conjecture}

\newtheorem{Ques}{Question}

\begin{document}
\title{Maharaja Nim \\ {\small Wythoff's Queen meets the Knight}}
\author{Urban Larsson and Johan W\"astlund,\\ {\small Mathematical Sciences},\\ {\small Chalmers University of Technology and University of G\"oteborg,}\\ {\small G\"oteborg, Sweden}\\ {\small urban.larsson@chalmers.se, wastlund@chalmers.se}}

\maketitle

\begin{abstract}
New combinatorial games are introduced, of which the most pertinent is Maharaja Nim. The rules extend those of the well-known impartial game of Wythoff Nim in which two players take turn in moving a single Queen of Chess on a large board, attempting to be the first to put her in the lower left corner. Here, in addition to the classical rules a player may also move the Queen as the Knight of Chess moves. We prove that the second player's winning positions of Maharaja Nim are close to the ones of Wythoff Nim, namely they are within a bounded distance to the lines with slope $\frac{\sqrt{5}+1}{2}$ and $\frac{\sqrt{5}-1}{2}$ respectively. For a close relative to Maharaja Nim, where the Knight's jumps are of the form $(2,3)$ and $(3,2)$ (rather than $(1,2)$ and $(2,1)$), we also demonstrate polynomial time complexity to the decision problem of the outcome of a given position. 
\end{abstract}

\section{Maharaja Nim}
We introduce a 2-player combinatorial game called 
\emph{Maharaja Nim}, an extension of the well-known game 
of Wythoff Nim \cite{Wy}. (The name `Maharaja' is taken from a variation 
of Chess, `The Maharaja and the Sepoys', \cite{Fa}.)   
Both these games are impartial, that is, the set of options 
are the same regardless of whose turn it is. For a background on impartial 
games see \cite{BeCoGu}. 

\begin{figure}[ht!]
\centering
\includegraphics[width=0.7\textwidth]{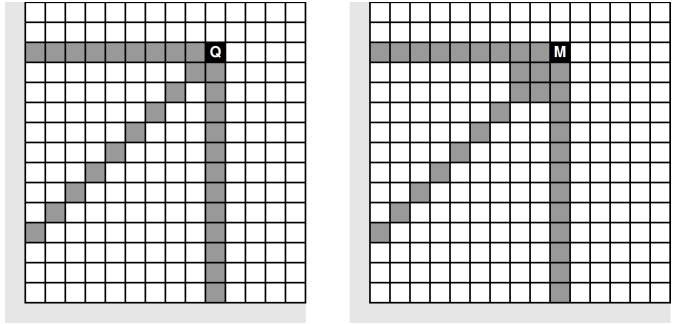}
\caption{The move options, from a given position, of Wythoff Nim and Maharaja Nim respectively.}\label{figure:1}
\end{figure}

Place a Queen (of Chess) on a given position of a large Chess board, with 
the position in the lower left corner labeled $(0,0)$. 
In the game of Wythoff Nim, here denoted by $\w$, the two players 
 move the Queen alternately 
as it moves in Chess, but with the restriction that, by moving, 
no coordinate increases, see Figure \ref{figure:1}. Only non-negative 
coordinates are allowed so that the first player 
to reach the position $(0, 0)$ wins. 

In Maharaja Nim, denoted by $\m$, the rules are as in Wythoff Nim, except 
that the Queen is exchanged for a `Maharaja', 
a piece which may move both as the Queen and the Knight of Chess,
again, provided by moving no coordinate increases. 
See Figure \ref{figure:1}.

\begin{figure}[ht]
\centering
\includegraphics[width=0.47\textwidth]{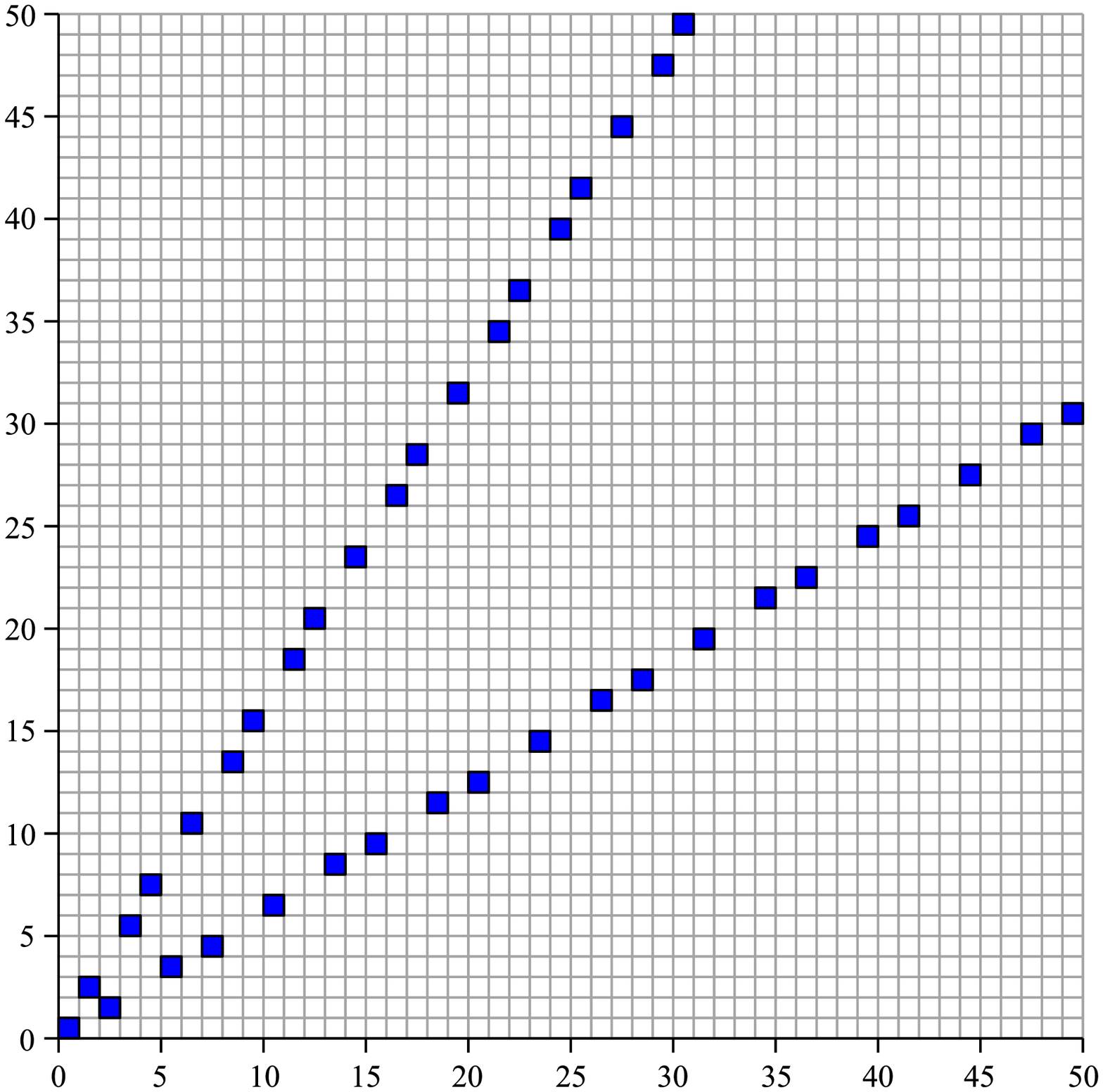}
\includegraphics[width=0.47\textwidth]{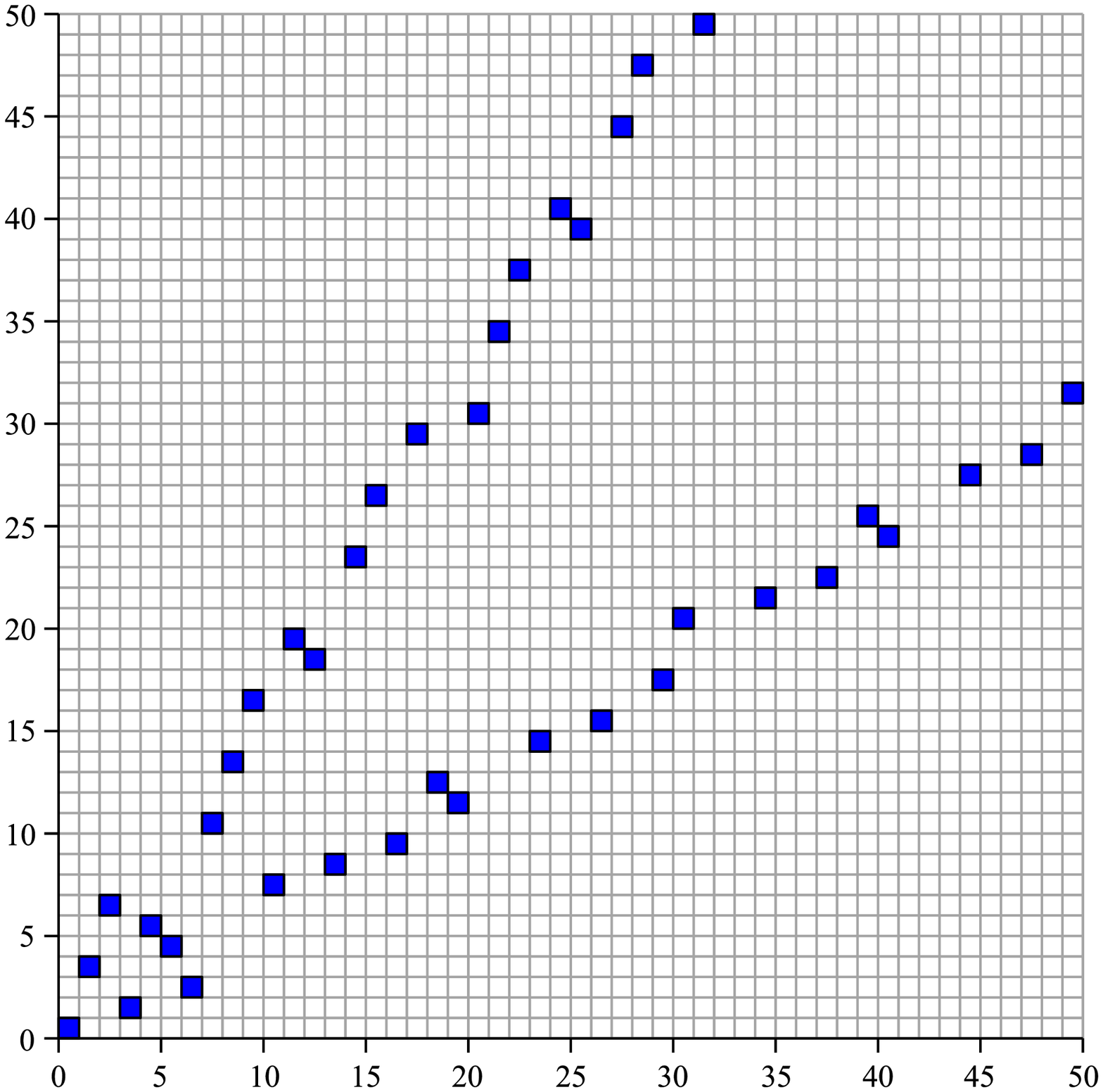}
\caption{The initial P-positions of Wythoff Nim and Maharaja Nim 
respectively.}\label{figure:2}
\end{figure}

We say that a position is P if the second player to move 
has a winning strategy, otherwise N. Also, let $\P_\m$ and $\P_\w$ denote 
the set of P-positions of Maharaja Nim and Wythoff Nim respectively. See Figure \ref{figure:2} for a computation of the initial P-positions of the respective games and the Appendix, Section \ref{A:1} for the corresponding code.

We let $\N$ denote the positive integers and $\M$ the non-negative integers. 
Let 
\begin{align*}
\phi = \frac{1+\sqrt{5}}{2}
\end{align*}
denote the golden ratio. The well-known winning strategy 
of Wythoff Nim \cite{Wy} is 
\begin{align}\label{ban}
\P_\w = \{(\lfloor \phi n \rfloor, \lfloor \phi^2 n \rfloor), 
(\lfloor \phi^2 n \rfloor, \lfloor \phi n \rfloor)\mid n\in \M\}. 
\end{align}
From this it follows that there is precisely one P-position 
of Wythoff Nim in each row and each column of the board (see also \cite{Be}).

The purpose of this paper is to explore the P-positions of Maharaja Nim 
and some related games. Clearly $(0, 0)$ is P. Another trivial observation 
is that, since the rules of game are symmetric, if $(x, y)$ is P then 
$(y, x)$ is P. It is also easy to see that there is at most one 
P-position in each row and each column 
(corresponding to the Rook-type moves). But, in fact, the same assertion 
as for Wythoff Nim holds:
\begin{Prop}\label{rowcolumn}
There is precisely one P-position of Maharaja Nim in each row and 
each column of $\M\times \M$.
\end{Prop}
\noindent{\bf Proof.} 
Since all Nim-type moves are allowed 
in Maharaja Nim, there is at most one P-position in each row and column 
of $\M\times \M$. This implies that there are at most $k$ P-positions 
strictly to the left of the $k^{th}$ column (row). Each such 
P-position is an option for at most 
three N-positions in column (row) $k$. This implies that 
there is a least position in column (row) 
$k$ which has only N-positions as options. By definition this 
position is P and so, since $k$ is an arbitrary index, 
the result follows. \hfill $\Box$\\

\begin{figure}[ht]
\centering
\includegraphics[width=0.45\textwidth]{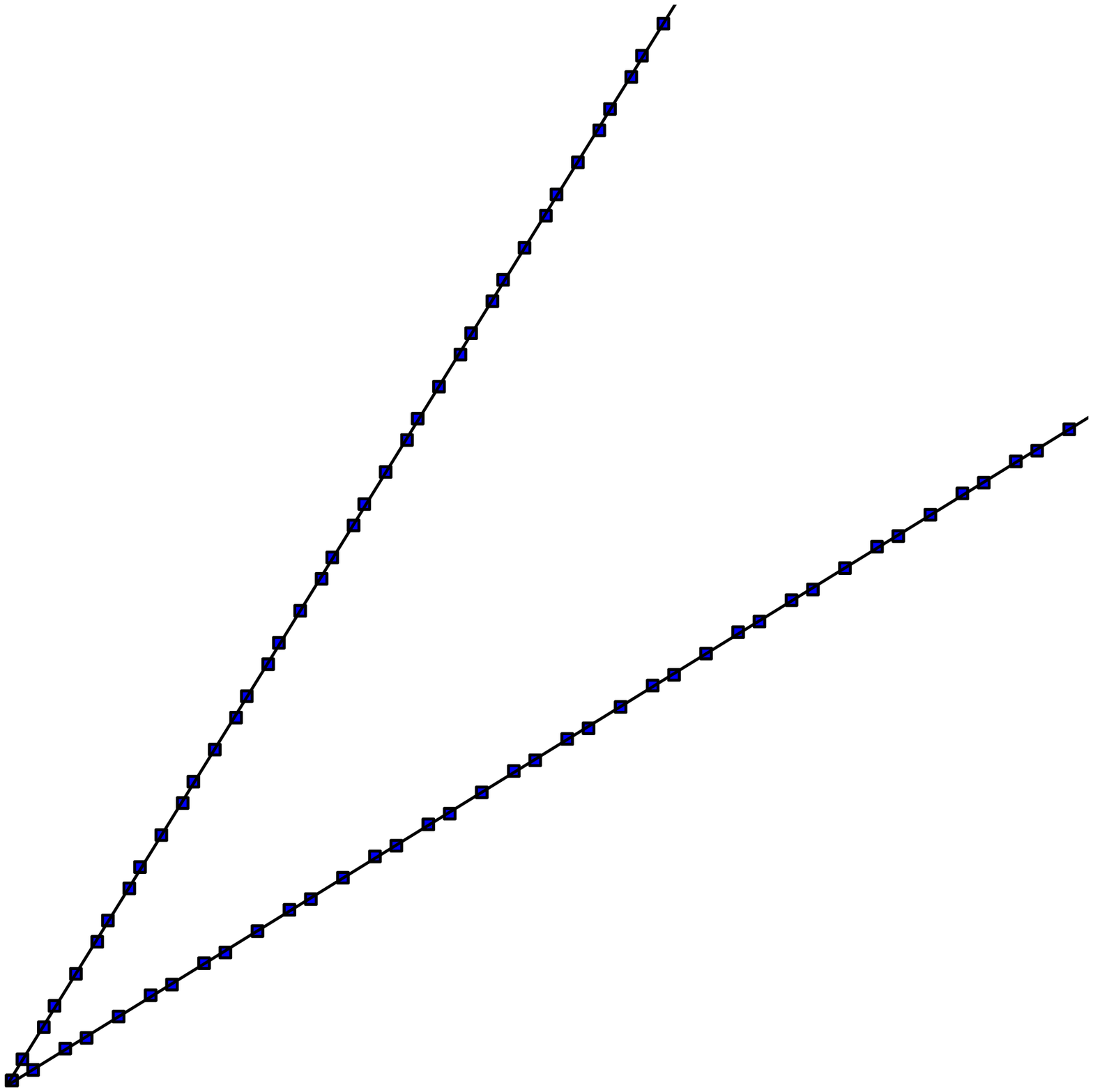}\includegraphics[width=0.45\textwidth]{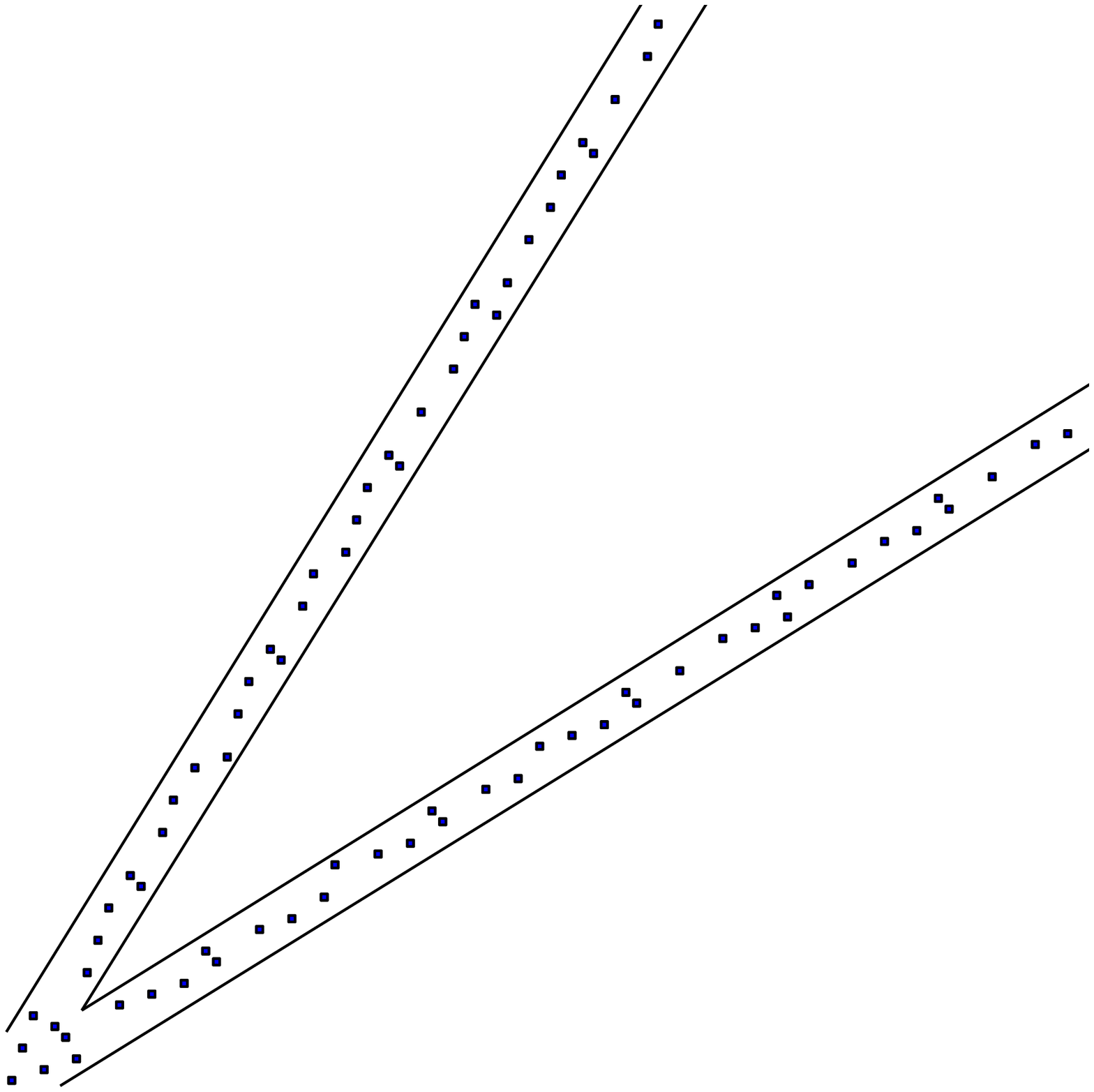}
\caption{To the left, the P-positions of Wythoff Nim lie `on' 
the lines $\phi x$ and $\phi^{-1}x$, $x\ge 0$. The figure to the right 
illustrates a main result of this paper, that the P-positions 
of Maharaja Nim are bounded below and above 
by the `bands', $\phi x + O(1)$ and $\phi^{-1} x + O(1)$}\label{figure:3}
\end{figure}

Another claim holds for both Wythoff Nim and Maharaja Nim. 
There is \emph{at most} one P-position on each \emph{diagonal} of the form 
\begin{align}\label{diag}
\{\{x, x + C\}\mid x\in \M\}, C\in \M,
\end{align}
(corresponding to the Bishop-type moves). But (\ref{ban}) 
readily gives that, for Wythoff Nim, there is \emph{precisely} one 
P-position on each such diagonal. Even more is true: If 
\begin{align}\label{abW}
\P_\w = \{(a_i,b_i),(b_i,a_i)\}, 
\end{align}
with $(a_i)$ increasing and for all $i$, $a_i\le b_i$, then for all $n$, 
\begin{align}\label{strengthening}
\{0, 1, \ldots , n\} = \{b_i - a_i\mid i\in \{0, 1,\ldots ,n\}\}.
\end{align}
As we will see in Section \ref{Sec:2}, a somewhat weaker, but crucial, 
property holds also for Maharaja Nim, but let us now state 
our main result (see also Figure \ref{figure:3}). We 
let $O(1)$ denote bounded functions on $\M$.
\begin{Thm}[Main Theorem]\label{mainthm} 
Each P-position of Maharaja Nim lies on one of the `bands' 
$\phi n + O(1)$ or $\phi^{-1} n + O(1)$, that is, if $(x, y)\in \P_\m$, 
with $y\ge x$, then $y - \phi x$ is $O(1)$.
\end{Thm}

We give the proof of this result in Section \ref{Sec:2}. It is 
quite satisfactory in one sense, but for the two gamesters trying 
to figure out how to quickly find safe positions, it does not 
quite suffice. The following question is left open.
\begin{Ques}\label{Ques:1}
Does Maharaja Nim's \emph{decision problem}, to determine whether a given 
position $(x, y)$, with input length $\log (xy)$, is P, 
have polynomial time complexity in $\log (xy)$?
\end{Ques}
In Section \ref{Section:5} we provide an affirmative answer of 
this question for a close relative 
of Maharaja Nim, namely the extension of Wythoff Nim where moves 
of type $(2,3)$ and $(3,2)$ are adjoined (but not $(1,2)$ or $(2,1)$). 
This result builds upon 
an analog result, of `approximately linear' P-positions, 
as that for Maharaja Nim in Theorem \ref{mainthm}. See also \cite{FrPe}, which was the inspiration for some results in this paper, 
although its methods do not seem to encompass the complexity of Maharaja Nim.

\subsection{Complementary sequences and a central lemma}
We say that two sequences of positive integers 
are \emph{complementary} if 
each positive integer is contained in precisely one of these sequences. In our setting this corresponds to Proposition \ref{rowcolumn} together with the claim before (\ref{diag}).
In \cite{FrPe} the authors proved the following result.

\begin{Prop}[Fraenkel, Peled]\label{AUthm}
Suppose $x$ and $y$ are complementary and increasing 
sequences of positive integers. Suppose further that 
there is a positive real constant, $\delta$, such that, for all $n$, 
\begin{align}\label{yxn}
y_n - x_n = \delta n + O(1). 
\end{align}
Then there are constants, $1<\alpha < 2 < \beta $, such that, for all $n$, 
\begin{align}\label{xn}
x_n - \alpha n = O(1)
\end{align}
and 
\begin{align}\label{yn}
y_n - \beta n = O(1).
\end{align}
\end{Prop}

As they have remarked (see also \cite{HeLa}), 
by simple density estimates one may decide the constants 
$\alpha$ and $\beta$ as functions of $\delta$. Namely, notice 
that (\ref{yxn}) and (\ref{xn}) together imply 
\begin{align}\label{betaalpha}
\beta = \alpha + \delta 
\end{align}
and, by complementarity, we must have  
\begin{align}\label{alphabeta}
\frac{1}{\alpha} + \frac{1}{\beta} = 1.
\end{align}
(Thus $\alpha$ and $\beta$ are algebraic numbers if and only if $\delta$ is.)
By this we get the relation
\begin{align}\label{deltalpha} 
\delta(1-\alpha) +\alpha = (\alpha -1)\alpha,
\end{align}
which will turn out to be useful in what will come next, namely we have 
found a short proof of an extension of their theorem---an extension which is  
easier to adapt to the circumstances of Maharaja Nim. Let us explain.

If we denote 
\begin{align}\label{abM}
\P_\m = \{(a_n, b_n),(b_n,a_n)\mid n\in \M\}, 
\end{align}
with $(a_n)$ increasing and for all $n, b_n\ge a_n$, then, for all $n, b_n$ 
is uniquely defined by the rules of M. At this point, one might want 
to observe that, if the $b$-sequence would have been increasing 
(by Figure \ref{figure:2} it is not) then 
Theorem \ref{mainthm} would follow from  Proposition \ref{AUthm}
if one could only establish the following claim: $b_n - a_n - n$ is $O(1)$.
Namely in (\ref{deltalpha}) $\delta = 1$ gives $\alpha = \phi$ 
in Proposition \ref{AUthm}. Now, interestingly enough, it turns out that Proposition 
\ref{AUthm} holds without the condition that the $y$-sequence is 
increasing, namely (\ref{yxn}) together with an increasing $x$-sequence 
suffices.

\begin{Lem}[Central Lemma]\label{centralthm}
Suppose $x$ and $y$ are complementary  
sequences of positive integers with $x$ increasing. Suppose further that 
there is a positive real constant, $\delta$, such that, for all $n$, 
\begin{align}\label{one}
y_n - x_n = \delta n + O(1). 
\end{align}
Then there are constants, $1 < \alpha < 2 < \beta $, such that, for all $n$, 
\begin{align}\label{xn2}
x_n - \alpha n = O(1)
\end{align}
and 
\begin{align}\label{yn2}
y_n - \beta n = O(1).
\end{align}
\end{Lem}

\noindent{\bf Proof.} 
We begin by demonstrating that, for all $n\in \N$, 
\begin{align}
x_{n+1} = x_n + O(1),\label{two}
\intertext{ and }
y_{n+1} = y_n + O(1).\label{three}
\end{align}
By (\ref{one}), for all $k,n\in \N$ we have that 
\begin{align}
y_{n+k} - y_n &= x_{n+k} + \delta(n+k) - x_n - \delta n + O(1),\notag\\
&= x_{n+k} - x_n + \delta k  + O(1).\label{yy}
\end{align}
Since for all $k, n \in \N$, 
$x_{n+k} - x_n \ge k$ and $\delta > 0$ this means that, for all $k, n \in \N$,
\begin{align}\label{yC}
y_{n+k}  \ge  y_n - C,
\end{align}
where $C$ is some universal positive constant (which may depend on $\delta$).
But, with $C$ as in (\ref{yC}), we can find another universal constant $\kappa = \kappa(C)\in \N$ such that, for all $n$, 
\begin{align}\label{2C}
y_{n+\kappa} - y_n \ge \kappa + 2C + 1.
\end{align} 
This follows since, in  (\ref{yy}), for any $C$, we can find $k=k(C)$ such that, for all $n$, $\delta k + O(1)> 2C$. Any such $k$ suffices as our $\kappa$. On the one hand there can be at most $\kappa - 1$ numbers from the $y$-sequence strictly between $y_n$ and $y_{n+\kappa}$   (with indexes strictly in-between $n$ and $n+\kappa$).
On the other hand the inequality (\ref{yC}) gives that there can be at 
most $C$ numbers from the $y$-sequence with index greater 
than $n + \kappa$ but less than $y_{n+\kappa}$. It also gives that there can be 
at most $C$ numbers with index less than $n$ but greater than $y_{n}$. 
Therefore, by complementarity and (\ref{2C}), there has to be a number
from the $x$-sequence in every interval of length $\kappa+2C+1$. Thus the jumps 
in the $x$-sequence are bounded, which is (\ref{two}). But 
then (\ref{three}) follows from (\ref{one}) and (\ref{two}) since 
\begin{align*}
y_{n+1} - y_n &= x_{n+1} + \delta(n + 1) - x_n - \delta n + O(1)\\
&=x_{n+1} + \delta - x_n + O(1)\\
&= O(1).
\end{align*}
By (\ref{three}) we may define $m$ as a function of $n$ with
\begin{align}\label{cons}
x_{n} = y_{m} + O(1).
\end{align}
(For example, one can take $m = m(n)$ the least number 
such that $x_n < y_m $. Then $y_m-x_n$ has to be bounded for otherwise $y_m-y_{m-1}$ is not bounded.) 
This has two consequences, of which the first one is 
\begin{align}\label{five}
x_n = n + m + O(1).
\end{align}
This follows since the numbers $1,2,\ldots , x_n$ are partitioned 
in $n$ numbers from the $x$-sequence, and the rest, by complementarity, $m + O(1)$ numbers from the $y$-sequence. 

The second consequence of (\ref{cons}) is that, by using (\ref{one}), 
\begin{align}\label{six}
x_n = x_m + \delta m + O(1).
\end{align}
If $\lim x_n/n$ and $\lim y_n/n$ exist then, clearly they must 
satisfy (\ref{betaalpha}) and (\ref{alphabeta}) with $\delta$ as 
in the lemma. Thus, using this definition of $\alpha = \alpha (\delta)$, 
for all $n$, denote  
\begin{align*}
\Delta_n := x_n - \alpha n.
\end{align*}
We want to use (\ref{five}) and (\ref{six}) to express 
$\Delta_n$ in terms of $\Delta_m$. 

Equation (\ref{six}) expresses $x_n$ in terms of $x_m$ and $m$. Therefore, 
we wish to combine (\ref{five}) and (\ref{six}) to express $n$ in terms 
of $x_m$ and $m$, that is, we wish to eliminate $x_n$ from (\ref{five}). 
If we plug in the expression (\ref{six}) for $x_n$ in (\ref{five}) and 
solve for $n$ we get  
\begin{align}\label{seven}
n = x_m + (\delta-1)m + O(1).
\end{align}
Combining (\ref{six}) and (\ref{seven}) gives
\begin{align}
\Delta_n &=  x_m + \delta m - \alpha(x_m + (\delta - 1)m) + O(1)\notag\\
&= (1-\alpha)x_m + (\delta (1-\alpha)+\alpha)m + O(1)\notag\\
&= (1-\alpha)\Delta_m + O(1),\label{eight}
\end{align}
where the last equality is by (\ref{deltalpha}).

Notice that, by (\ref{six}), for sufficiently large $n$ we have that $m < n$. 
Hence we may use strong induction and by (\ref{eight}) conclude 
that $\Delta_n$ is $O(1)$ which is (\ref{xn2}). 
Then (\ref{yn2}) follows from (\ref{one}).
\hfill $\Box$\\

\section{Perfect sectors, a dictionary and the proof 
of Theorem \ref{mainthm}}\label{Sec:2}
This whole section is devoted to the proof of Theorem \ref{mainthm}. We begin 
by proving that there is precisely one P-position of Maharaja Nim 
on each diagonal of the form in (\ref{diag}). 
Then we explain how the proof of this result leads to the second part 
of the theorem, the bounding of the P-positions within the 
desired `bands' (Figure \ref{figure:3}).

A position, say $(x,y)$, is an \emph{upper} position if it is strictly 
above the \emph{main diagonal}, that is if $y > x$. 
Otherwise it is \emph{lower}.

\begin{figure}[ht]
\centering
\includegraphics[width=0.4\textwidth]{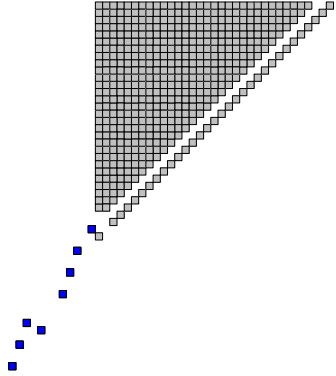}
\caption{All upper positions from which a player 
can move to an upper P-position are erased. (The `sector' continues 
above the figure.) However, the `sector' is not perfect.}
\label{figure:5}
\end{figure}

\begin{figure}[ht]
\centering
\includegraphics[width=0.4\textwidth]{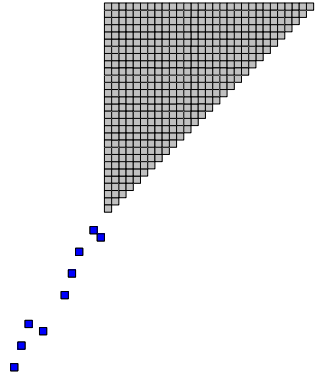}
\caption{(Step 1) A perfect sector together 
with the corresponding initial P-positions.}
\label{figure:6}
\end{figure}

We call a \emph{$(C,X)$-perfect sector}, or simply a \emph{perfect sector},
all positions strictly above some diagonal of the form in (\ref{diag}) 
and strictly to the right of column $X$.
Suppose that we have computed all P-positions in the 
columns $1,2, \ldots , a_{n-1}$ and that, when we 
erase each upper position from which a player can move to an upper P-position 
(Figures \ref{figure:5} and \ref{figure:6}), then the remaining 
upper positions strictly to the right of $a_{n-1}$ constitute 
an $(n-1,a_{n-1})$-perfect sector (Figure \ref{figure:6}). Then we say that $a_{n-1}$ is \emph{perfect} and, in fact, it is easy to see that also property (\ref{strengthening}) holds for all such $n$. On the other hand, we will see that the converse statement holds if and only if for any such $n$,
\begin{align}\label{nth}
b_{n} - a_{n} = n,
\end{align}
given that the lower P-positions do not interfere. It is crucial to our approach that the first implication can be made stronger to also include (\ref{nth}). 
\begin{Lem}\label{LemmaPerfect}
Let $n\in \N$ be sufficiently large so that Knight type moves from lower P-positions do not affect the coordinates of upper P-positions and define $(a_i)$ and $(b_i)$ as in (\ref{abM}). Suppose also that
\begin{align}\label{12n}
\{0, 1, \ldots , n-1\} = \{b_i - a_i\mid 0\le i < n \} 
\end{align}
holds. Then (\ref{nth}) must hold if and only if $a_{n-1}$ is perfect.
\end{Lem}

\noindent{\bf Proof.} Suppose that (\ref{nth}) does not hold. Then clearly $b_n-a_n>n$. This must be due to a Knight type move from an upper P-position from $(a_{n-1},b_{n-1})$, that is to position $(a_{n-1}+1,b_{n-1}+2)$. Hence $a_{n-1}$ is not perfect. For the other direction, whenever there is no $i<n$ such that $b_i=a_{n-1}+1$, so that $a_n=a_{n-1}+1$ excludes a Knight type move as in the previous paragraph and hence assures a perfect sector.
\hfill $\Box$\\

\subsection{Constructing Maharaja Nim's bit-string}
We study a bit-string, a sequence of `0's and `1's, where the 
$i^{th}$ bit equals `0' if and only if there is an upper 
P-position of Maharaja Nim in column $i$. 
By Proposition \ref{rowcolumn}, if there is no upper P-position in 
column $i$, there is a lower ditto (the $i^{th}$ bit equals $1$).

\begin{figure}[ht]
\centering
\includegraphics[width=0.45\textwidth]{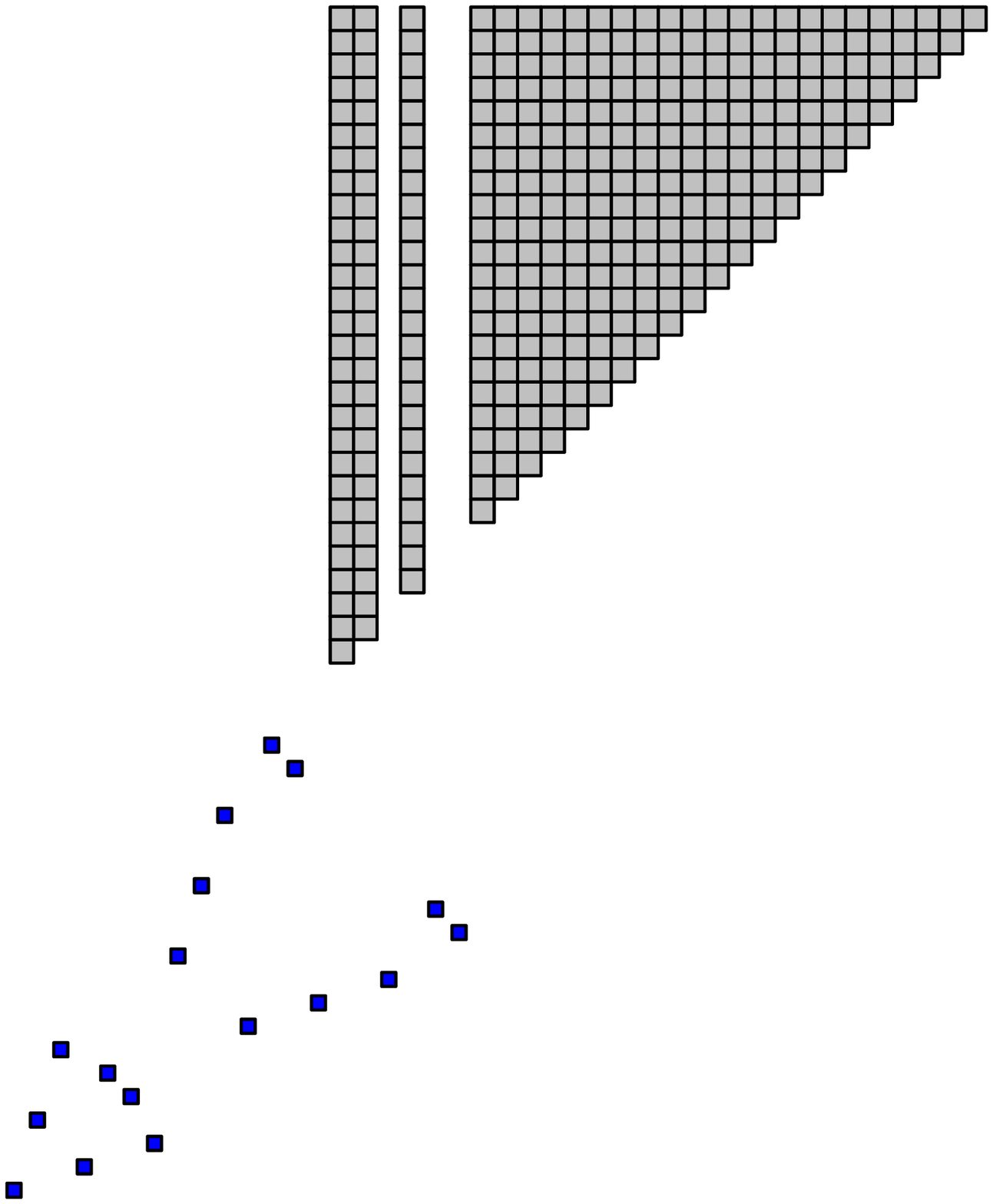}
\caption{
(Step 2) Each column in the perfect sector which
corresponds to a lower P-position (a `1' in the bit-string) 
has been erased.}
\label{figure:7}
\end{figure}

Suppose that $a_x=n$ is perfect. Then, by symmetry we know some 
lower P-positions in columns to the right of $n$. The next step 
is to erase each column in this perfect sector which has a 
lower P-position, a `1' in the bit-string (see Figure \ref{figure:7}) and 
to, recursively in the non-erased part of the perfect sector, compute new 
upper P-positions. We do this until we reach the next perfect 
sector (for the moment assume that this will happen) at say column 
$n + m$, $m > 0$. Thus, using this notation, we 
may define a word of length $m$, containing the information of whether 
the P-position in column $i\in \{n, n+1, \ldots , n+m-1\}$ is below or 
above the main diagonal. 

At this point we adjoin this \emph{word} together with its unique 
\emph{translate} to Maharaja Nim's \emph{dictionary}.
The translate is obtained accordingly: For each P-position in the columns 
$n$ to $n + m - 1$ define the $i^{th}$ bit in the translation as a `1' 
if and only if row $k + i$ has an upper P-position and where $k$ is the 
largest row index strictly below the perfect sector. See 
also Figure \ref{figure:8} and the next section for examples.
Then the translate has length $m + l$, where $l$ denotes the number of  
`0's in the word.

We then concatenate the translate at the end of the existing bit-string. 
In this way, provided a next perfect sector will be detected, 
the bit-string will always grow faster than we read from it.
However, there is no immediate guarantee that we will be able to repeat 
the procedure---that the next word exists---or 
for that matter that the size of 
the dictionary will be finite, so that the process may be described by a 
finite system of words and translates. But, in the coming, 
we aim to prove that, in fact, the next perfect sector will always  
(in the sense outlined above) be detected 
within a `period' of at most 7 P-positions, that is `0's in the bit-string. 
As we will see, a complete dictionary needs only (between 9 and) 
14 translations.

Let us describe a bit more in detail how the first part of Maharaja's 
bit-string is constructed. 

\subsection{A detailed example}\label{Section:2.2} 
Initially there is some interference which does not allow a recursive 
definition of words and translates, see Figure \ref{figure:2}.
The first perfect sector beyond the origin is attained 
when the four first P-positions strictly above the main diagonal has 
been computed. This happens to the right of column 8. To the right of 
column 12 a new perfect sector is detected. Thus the first word 
(left hand side entry) in the 
dictionary will be `00100', corresponding to the P-positions 
$(8,13)$, $(9,16)$, $(10,7)$, $(11,19)$ and $(12, 18)$. (Here there is 
no interference since the $y$-coordinate of the first P-position 
is greater than the $x$-coordinate of the last P-position.) 
Let us verify that this word translates to `100101100'.
Notice that the first `1'-bit means that the P-position $(8, 13)$ is 
to the left of the main diagonal---by symmetry this corresponds to a 
lower P-position in column 13. The second bit is `0'. This means that 
the next upper P-position is in column 14. Then, by rules of game, 
it has to be at least in row $16$, which indeed will be attained, so 
that the next P-position will be $(9, 16)$. By the rules of game, the 
rows 14 and 15 cannot have P-positions to the left of the main diagonal, 
so that a prefix is `1001'. Continuing up to the last P-position of 
this translate, $(12, 18)$, extends the prefix to `1001011'. The next 
upper P-position will be in at least row $22$ since the least 
unused diagonal is $22 - 13 = 9$. 
After this a new perfect sector will start. This gives the two last `0's in 
the translate, `100101100', which may now be concatenated at the end 
of the first part of the bit-string, `00100', so that 
the new bit-string becomes `00100\underline{1}00101100'. 

\begin{figure}[ht!]
\centering
\includegraphics[width=0.4\textwidth]{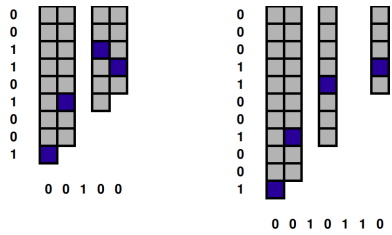} 
\caption{To the left, the unique (upper) P-positions of Maharaja Nim 
in the columns 8 to 12 are computed. The corresponding translation is $00100\rightarrow 100101100$. To the right are the P-positions in the columns 14 to 20 together with the translate $0010110\rightarrow 10010011000$. (Here we 
have omitted column 13 with its translation $1\rightarrow 0$.) 
See also Figure \ref{figure:2} and Section~\ref{Section:2.2}.}
\label{figure:8}
\end{figure}

In column 13 there is a lower P-position (corresponding to the $6^{th}$ 
bit in the string), which gives a new perfect sector by default,
that is, the next left hand side word is  `1'. This corresponds 
to that the first column in a perfect sector is erased and we get a 
new perfect sector without adding any upper P-position. 
By the property of a perfect sector, there can be 
no P-position to the left of the main diagonal in row 22, so 
the translate of the word `1' must be `0'. A concatenation of this `0' 
at the end of the existing string gives `001001\underline{0}01011000'. 
As we continue to read from the `0' in the seventh position it turns out that, 
this time, we need to read `0010110' (Figure \ref{figure:8} to the right) 
to obtain a new perfect sector and also 
that this word translates to `10010011000'. Again, concatenating 
this translate at the end of the 
existing string gives `0010010010110\underline{0}010010011000', and so on.\\\\

\subsection{Maharaja Nim's dictionary}
The dictionary of $M$ is
\begin{align}
1 &\rightarrow 0\label{dic1}\\
01 &\rightarrow 100\\
00100 &\rightarrow 100101100\\
00110 &\rightarrow 10010100\\
000100 &\rightarrow 10010110100\\
001110 &\rightarrow 100100100\\
0010110 &\rightarrow 10010011000\\
00000100 &\rightarrow 100101100111000\\
000010010 &\rightarrow 1001001111000100\label{dic9}\\
0000000 &\rightarrow 10010110110100\\
0010100 &\rightarrow 100100110100\\
0011110 &\rightarrow 1001000100\\
00000010 &\rightarrow 100101101100100\\
00001000 &\rightarrow 100100111100100.
\end{align} 

By computer simulations we have verified that each one of the words (\ref{dic1}) to (\ref{dic9}) does appear in Maharaja Nim's bit-string. We have included the code the Appendix, Section \ref{A.2}. By our method of proof, we have found no way to exclude the latter five, but a guess is that they do not appear. At least they do not appear among the first 20000 bits of the bit-string. The following result gives the first part of the theorem.\\

\begin{Lem}[Completeness Lemma]\label{Lemma}
When we read from Maharaja Nim's bit-string each prefix is 
contained in our extended dictionary of (left hand side) words 
of Maharaja Nim. 
\end{Lem}
\newpage
\noindent{\bf Proof.}
Let us present a list in lexicographic order of all words in 
our extended dictionary together with the words we need to exclude:
\begin{align*} 
0000000 &\rightarrow 10010110110100\\
00000010 &\rightarrow 100101101100100\\
00000011 &\text{ 'to exclude' (a)}\\
00000100 &\rightarrow 100101100111000\\
00000101 &\text{ 'to exclude' (b)}\\
0000011 &\text{ 'to exclude (c)'}\\
00001000 &\rightarrow 100100111100100\\
000010010 &\rightarrow 1001001111000100\\
000010011 &\text{ 'to exclude' (d)}\\
0000101 &\text{ 'to exclude' (e)}\\
000011 &\text{ 'to exclude' (f)}\\
000100 &\rightarrow 10010110100\\
000101 &\text{ 'to exclude' (g)}\\
00011 &\text{ 'to exclude' (h)}\\
00100 &\rightarrow 100101100\\
0010100 &\rightarrow 100100110100\\
0010101 &\text{ 'to exclude' (i)}\\
0010110 &\rightarrow 10010011000\\
0010111 &\text{ 'to exclude' (j)}\\
00110 &\rightarrow 10010100\\
001110 &\rightarrow 100100100\\
0011110 &\rightarrow 1001000100\\
0011111 &\text{ 'to exclude' (k)}\\
01 &\rightarrow 100\\
1 &\rightarrow 0
\end{align*}
This list is `complete' in the sense that any bit-string has precisely 
one of the words on the left hand side as a prefix. 
This motivates why it suffices to exclude the words `to exclude'. For example 
(a) needs to be excluded since the only word in our list beginning 
with `0000001' ends with a `0'. Neither can we translate words beginning 
with `000001's and ending with `01' or `1'. This motivates 
why we need to exclude (b) and (c). All left hand side 
words in our dictionary beginning with 4 `0's continues with 
100, which motivates that (e) and (f) need to be excluded, and so on. 
We move on to verify that the strings (a) to (k) are not contained
in the bit-string. 

No translate can contain more than three consecutive `0's. 
To get a longer string one has to finish off one translate and start a new. 
The only translate which starts with `0' is `0'. Thus, when a sequence of 
four or more `0's is interrupted it means that a new translate has begun. 
But all translates that begin with a `1' begins with `100'. Thus, a 
sequence of four or more `0's cannot be followed by `11' or `101'. This 
gives that the exclusion of the words (a),(b), (c), (e) and (f) is correct. 

Clearly, the string `100' in (d) has to be the prefix of some translate. Since 
the next two bits are `11', by the dictionary, this translate has to 
be `100'. But then the next translate has the prefix `11', 
which is impossible.

For the exclusion of (g) and (h) notice that the only strings of 
three consecutive `0's 
that exist within a translate is either at the end or is followed 
by the string `100'. Therefore, a string of three `0's cannot be followed 
by `11' or `101'. 

For (i), notice that the sub-string `101010' is not contained 
in any translate. If it were, it needed to be either at 
the beginning of a translate, which is impossible (since all of them 
except `0' begin with `100') or be split between two. 
The latter is impossible since all translates except `0' ends with `00'. 
In analogy to this, also (j) must be excluded and similarly for 
(k) since no translate contains 5 consecutive 
`1's and all translates ends in a `0', but starts with 
either `0' or `10'.\hfill $\Box$\\

Since the left hand side words have at most 7 `0's we adjoin at 
most 6 P-positions in a sequence with $b_n - a_n$ distinct from $n$. Namely, 
by Lemma \ref{LemmaPerfect}, when we start a 
new perfect sector we know that the next P-position will satisfy 
$b_n - a_n = n.$ The number of bits in a translate is bounded (by $16$) 
so that $b_n$ can never deviate more than a bounded number of 
positions from $a_n + n$. Hence, by Proposition \ref{rowcolumn}, 
the conditions of Lemma \ref{centralthm} are satisfied with 
the $a$-sequence as $x$, the $b$-sequence as $y$ and $\delta = 1$. Thus, 
 $b_n- a_n - n$ is $O(n)$ (as discussed in the paragraph 
before Lemma \ref{centralthm}) this concludes the proof 
of Theorem \ref{mainthm}. By inspecting the dictionary one can see that, in fact, for all $n$, $-4 \le b_n-a_n -n\le 3$.

\section{Dictionary processes and undecidability}
Let us briefly discuss a problem related to the method used in this paper.
Given a dictionary (of binary words and translations) and a starting string, 
will the translation process of the bit-string `terminate' or not?

More precisely, let us assume that we have a finite 
list of words $A = \{A_1, A_2, \ldots , A_m\}$ with 
translates $B_1, B_2, \ldots , B_m$ respectively, 
each word being a string of `0's and 
`1's, and where we assume that none of the words in $A$ is a prefix 
of another. (The latter is a convenient, but not necessary, condition 
as we will explore further in Section \ref{Section:5}.) Namely, as 
the read head reads from the bit-string, a natural generalization of 
a prefix free dictionary is to translate precisely the longest word 
containing a certain prefix.
 
Take any string $S$ as a starting string (for example $A_1$ but it could 
be an arbitrary string, not necessarily in the list). A `read head' `$\_$' 
starts to read $S$ from left to the right and as soon as it finds a 
string $A_i$ in $A$ it stops, sends a signal to a printer at the other end 
which concatenates the translation $B_i$ at the end of $S$. Then the read 
head continues to read from where it ended until it finds the next word 
in $A$, its translation being concatenated at the end, and so on.

If the read head gets to the end of the string without finding a word in 
the list $A$, the process stops with the current string as `output'. 
Otherwise, the process continues and gives as output an infinite string.

It follows from E. Post's tag productions \cite{Mi, Po} 
that it is algorithmically undecidable whether our 
`dictionary processes' stop or not. We give a proof in the Appendix, part B.

\section{Approximate linearity, converging dictionaries 
and polynomial time complexity}\label{Sec:4}
\begin{figure}[ht!]
\centering
\includegraphics[width=0.17\textwidth]{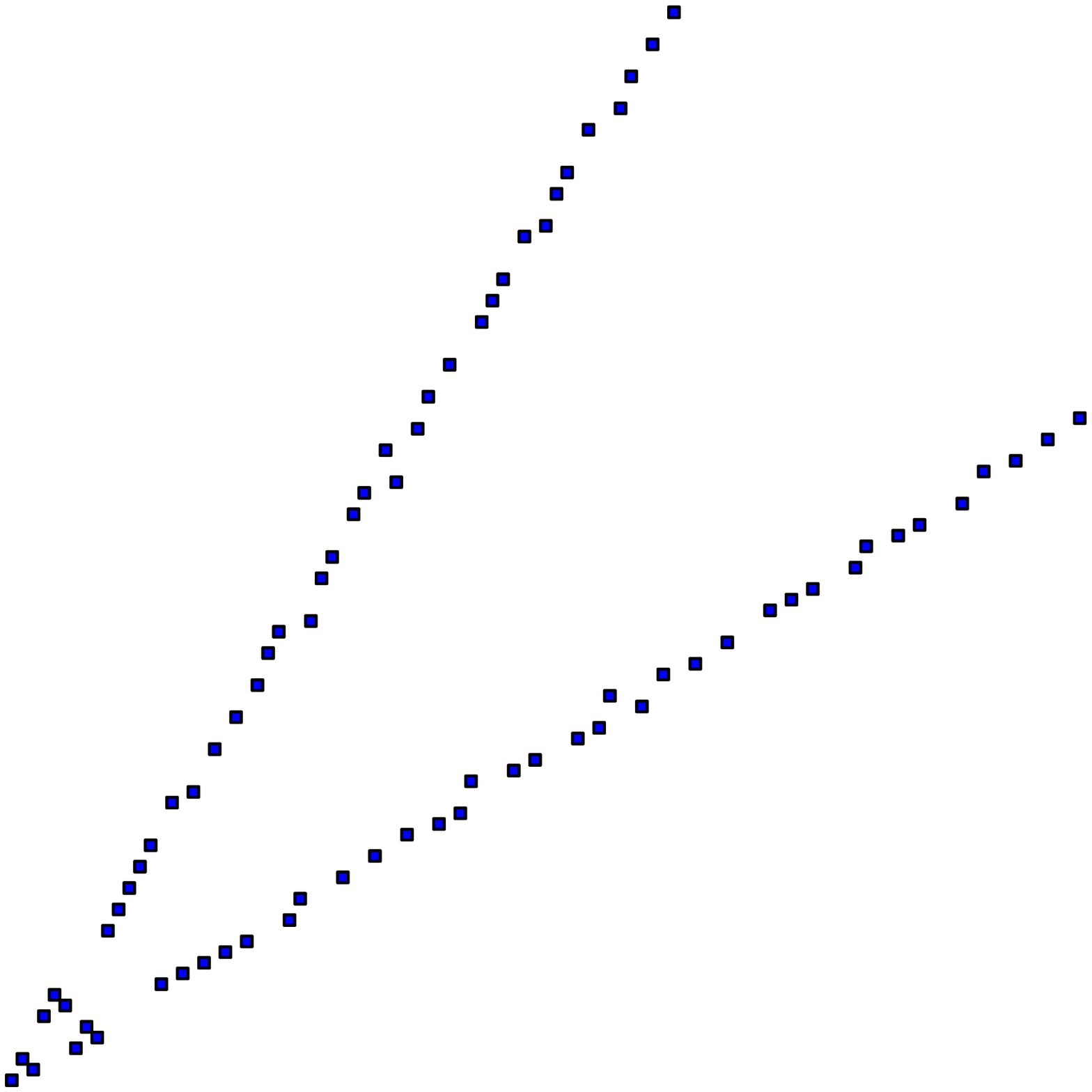}\includegraphics[width=0.17\textwidth]{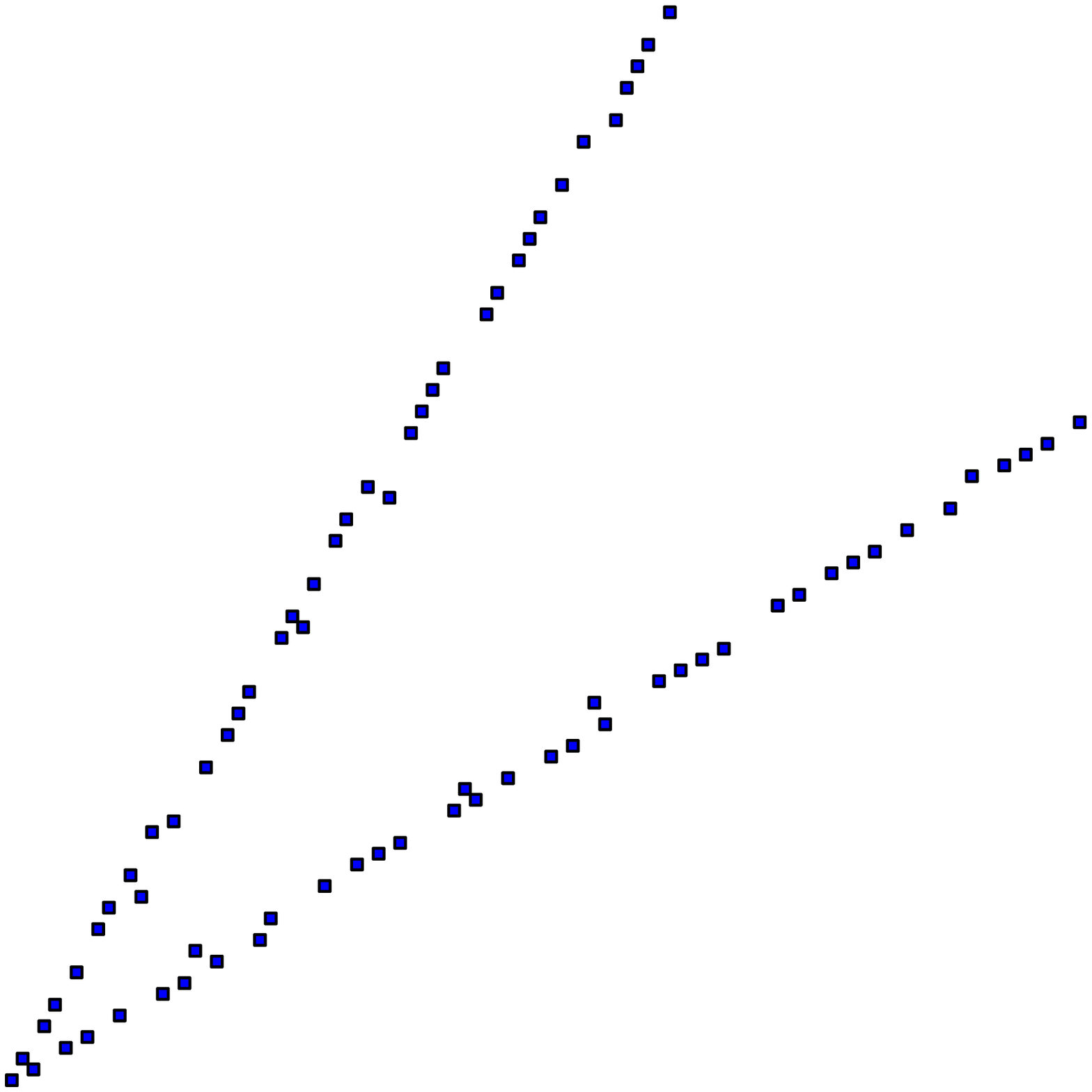}\includegraphics[width=0.17\textwidth]{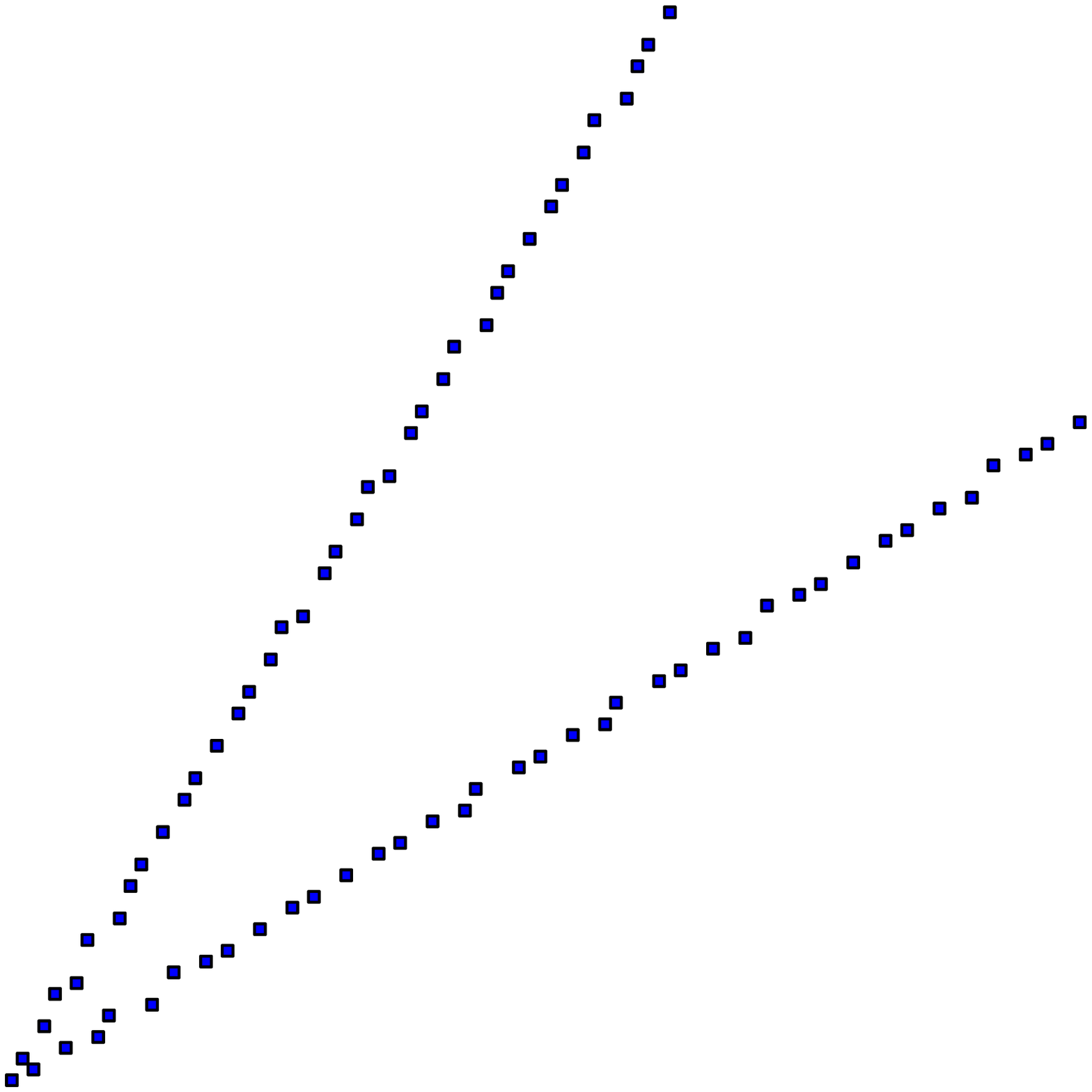}\includegraphics[width=0.17\textwidth]{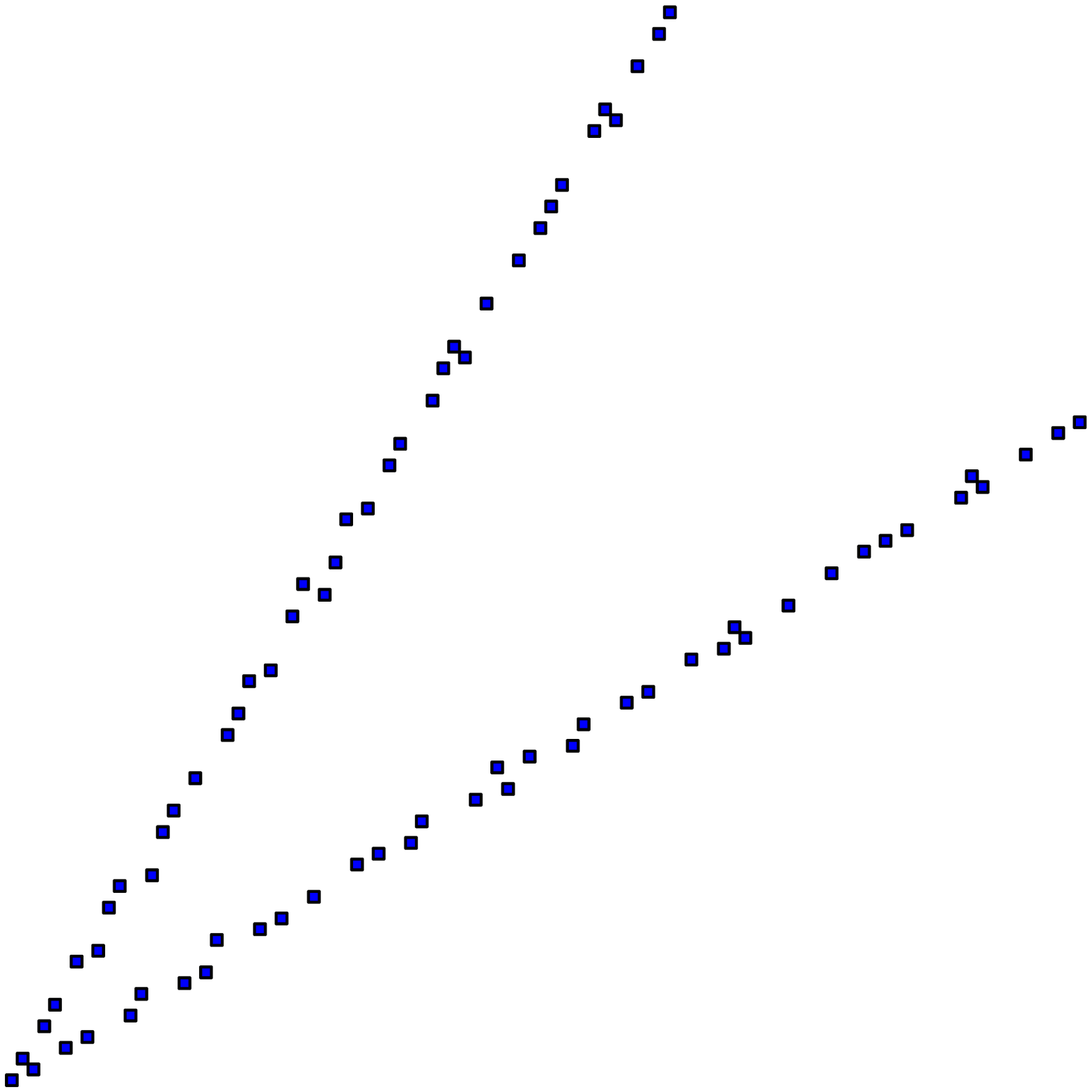}\includegraphics[width=0.17\textwidth]{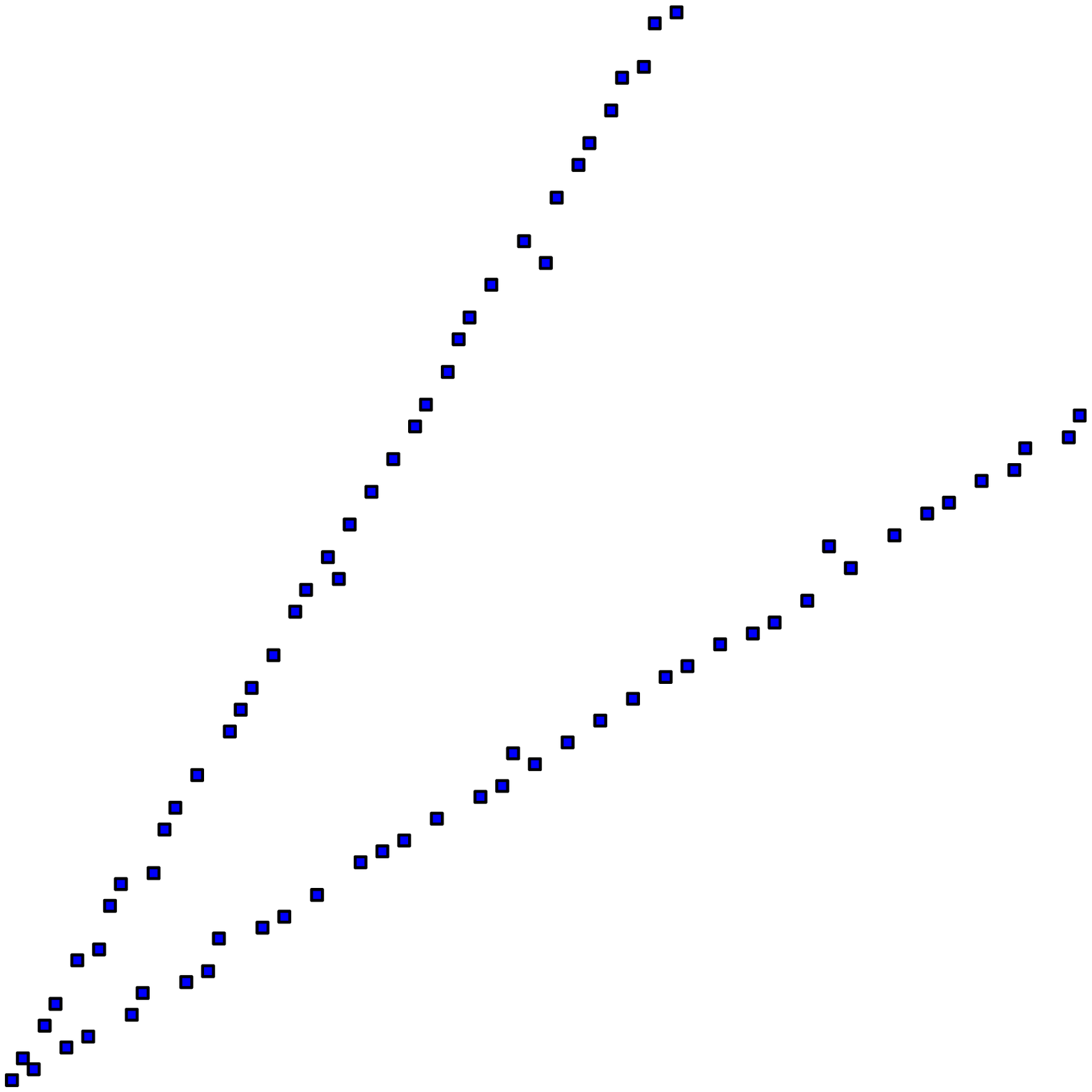}\includegraphics[width=0.17\textwidth]{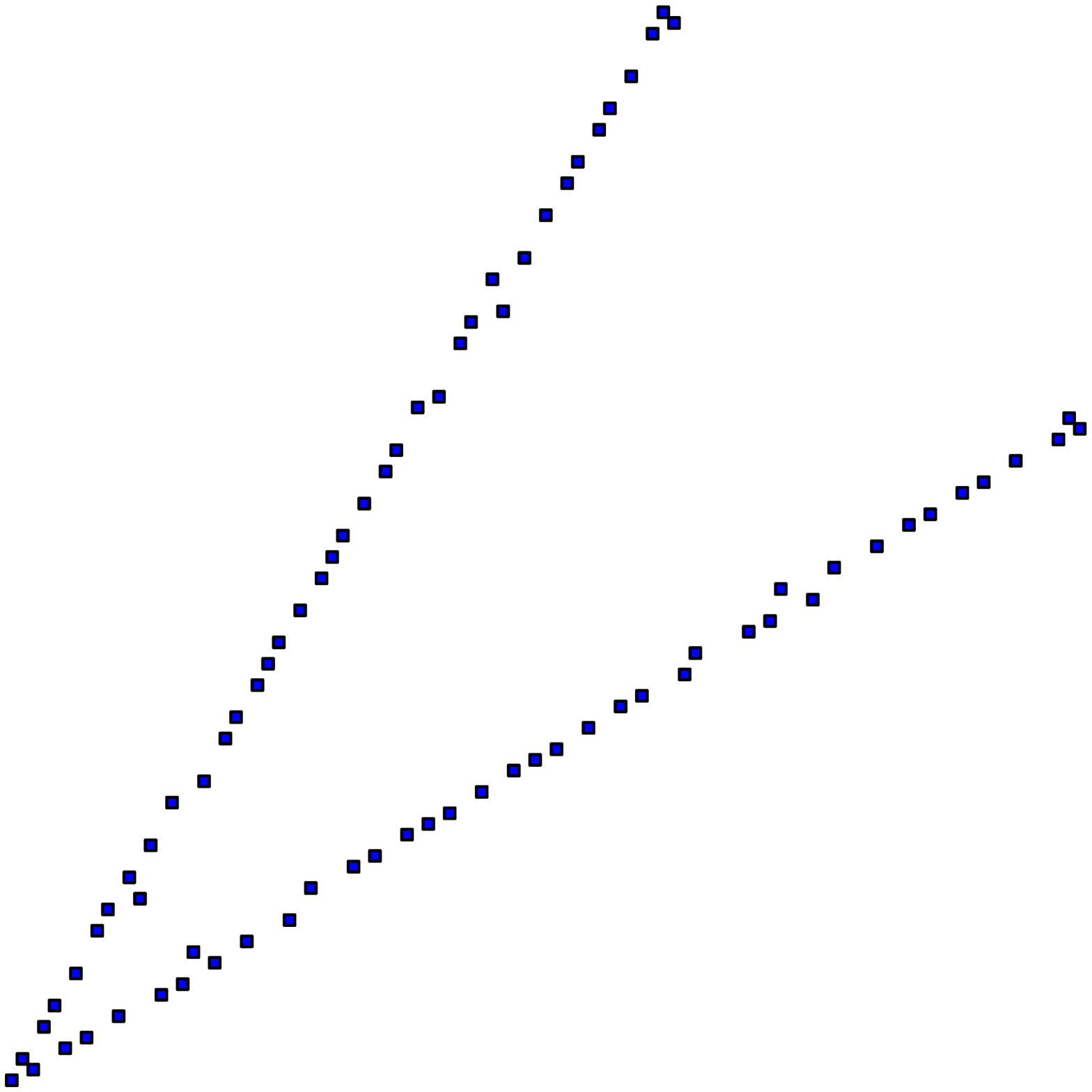}
\caption{The initial P-positions (the coordinates are less than 100) 
of $(k,l)$M for $(k,l) = (3,5),(4,6),(4,7),(5,8),(6,10)$ and $(7,11)$ 
respectively. In support of Conjecture \ref{conjecture:1}, the ratios of 
the respective coordinates seem to closely approximate $\phi$ or $1/\phi$.
(For $(2,3)$M, see Section \ref{Section:5}.)}
\label{figure:10}
\end{figure}

\begin{figure}[ht!]
\centering
\includegraphics[width=0.25\textwidth]{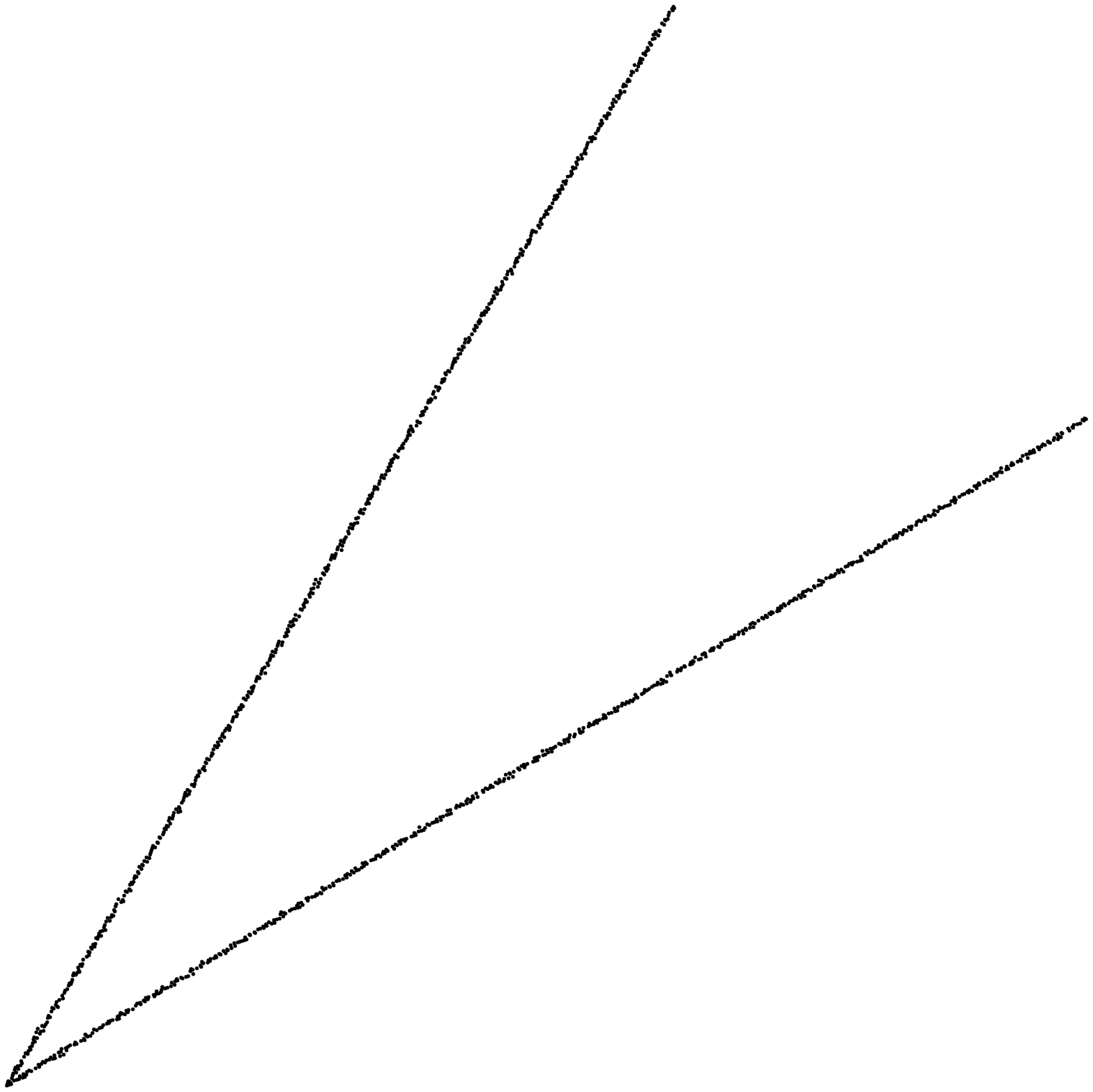}\includegraphics[width=0.25\textwidth]{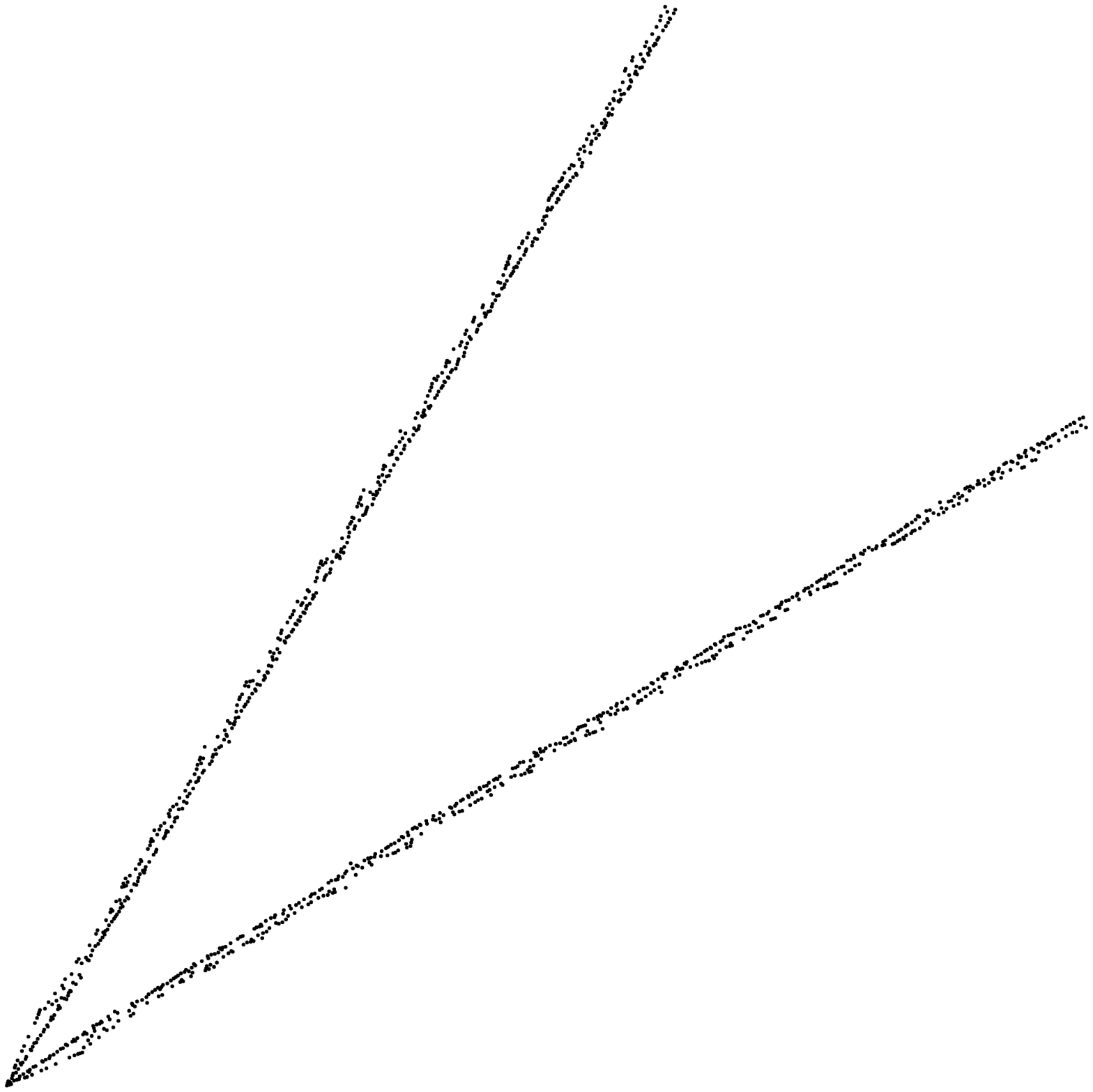}\includegraphics[width=0.25\textwidth]{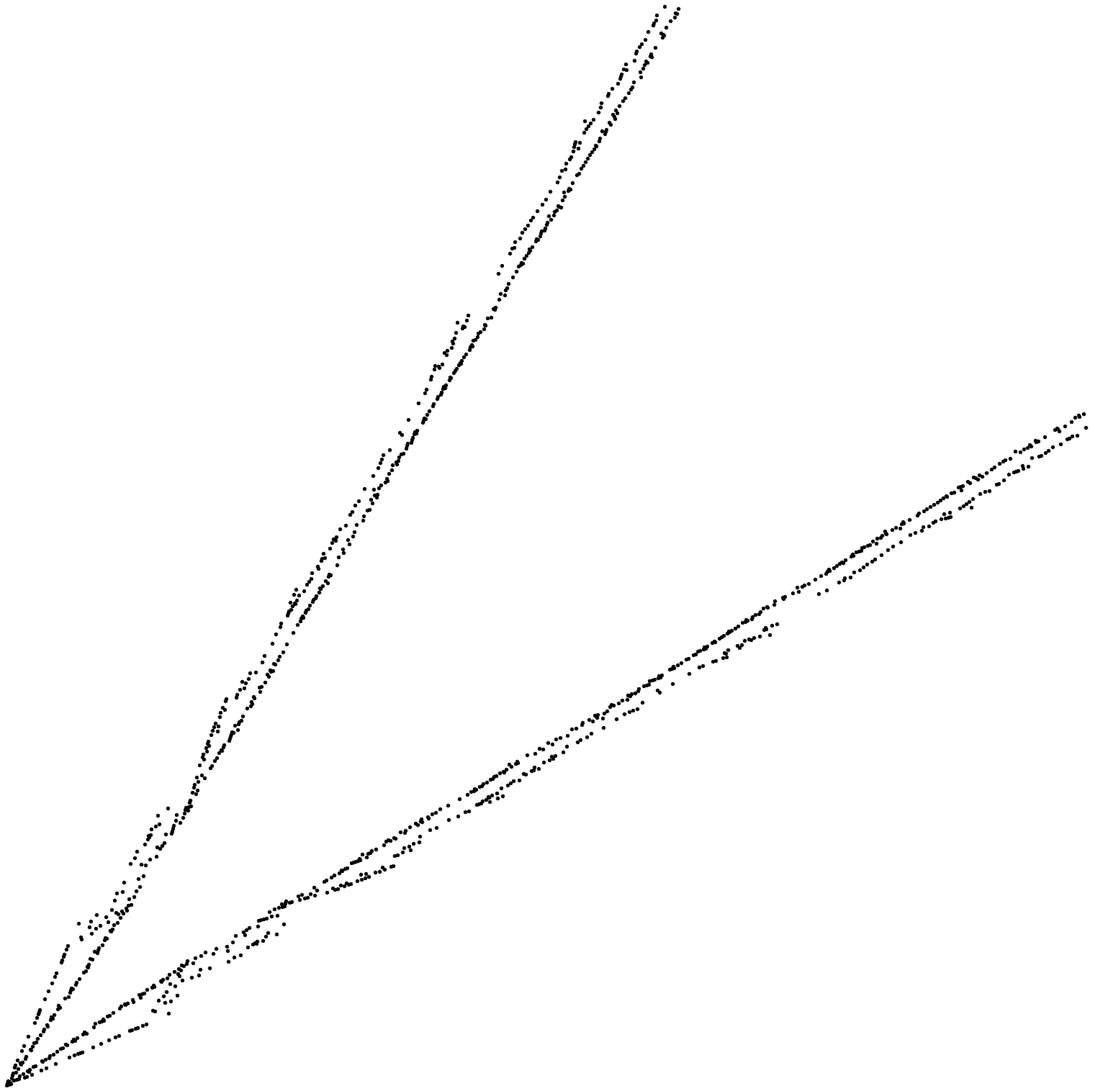}\includegraphics[width=0.25\textwidth]{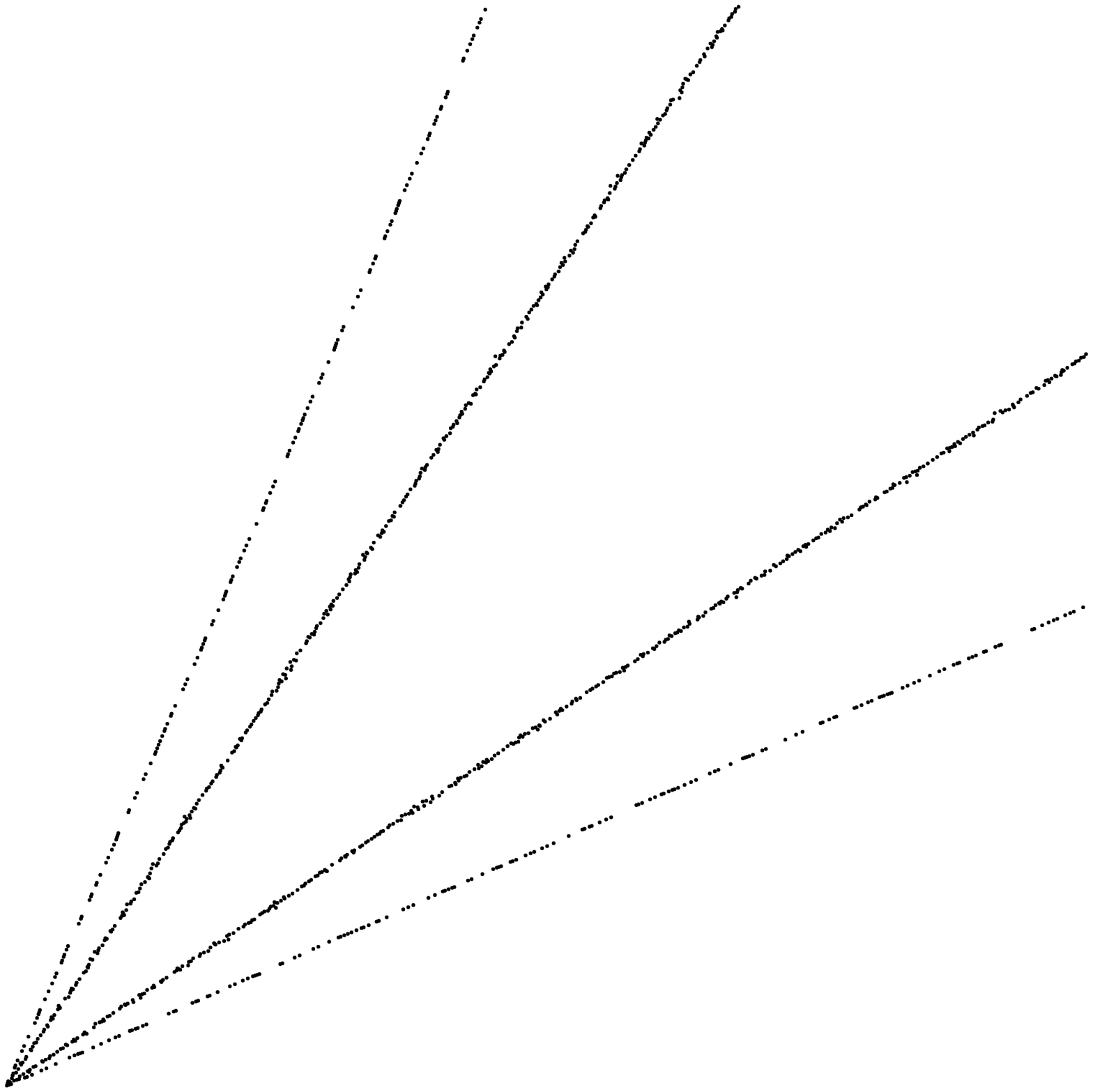}
\caption{The initial P-positions (the coordinates are less than 1500) 
of four extensions of Maharaja Nim where the adjoined moves are 
$\{(t,2t),(2t,t)\}$ where $t\in \{1, 2, \ldots, 10\}$, $\{1, 2, \ldots, 50\}$, 
$\{1,2,\ldots, 100\}$ and $\N$ respectively. That is the three first games 
have a finite number of moves adjoined to Wythoff Nim but the last one 
has infinitely many. Notice the seemingly emerging `bounded split' of the (upper) P-positions in the middle two figures, the ratio of the coordinates still seem to be within a bounded distance of $\phi$, but in the last figure, where an infinite number of moves are adjoined the convergence to $\phi$ is destroyed, a fact which is proved in \cite{La12}, and an 'unbounded split' (as in the rightmost figure) is established, which was recently proved in \cite{La}).}
\label{figure:11}
\end{figure}

There are infinitely many relatives to Maharaja Nim of the form `adjoin 
a finite set of moves to Wythoff Nim'. It is easy to see that Proposition \ref{rowcolumn} and (\ref{diag}) hold also for these type of games. 
For any given such generalization, is it possible to determine the greatest 
departure from $n$ for  $b_n - a_n$? For example see the games in 
Figure \ref{figure:10} and \ref{figure:11}.
Even simpler, is it decidable, whether there is a P-position above some 
straight line? More precisely: 
\begin{Ques}\label{Ques:2}
Given the moves of Wythoff Nim together with 
some finite list of moves, that is 
ordered pairs of integers (in Maharaja Nim the list is $\{(1, 2), (2, 1)\}$) 
and a linear inequality in two variables $x$ and $y$, is it 
decidable whether there is a P-position in the game which satisfies 
the inequality? 
\end{Ques}

On the one hand it is not even clear if 
a `generalized Maharaja Nim' has a finite dictionary in the sense 
of Section \ref{Sec:2}. On the other hand the solution of a similar 
game may or may not depend on the possible outcome of a dictionary process 
as in Section \ref{Sec:2}. In fact, in Section \ref{Section:5} we prove that 
a related dictionary process is successful in giving a polynomial time 
algorithm for the decision problem of whether a certain position is P. 
Therefore, let us look into some questions regarding some close 
relatives of Maharaja Nim.

To begin with, one might want to pay special attention to the 
family of extensions of Wythoff Nim, where the 
adjoined moves are of the form $(k, l)$ and $(l, k)$, 
$k, l\in \N$, $k < l$.  We call a game in this family $(k,l)$-Maharaja Nim, 
$(k,l)$-M. (Another problem is indicated in Figure \ref{figure:11} and its discussion.) The P-positions are distinct from those of Wythoff Nim, see \cite{La}, if and only if $(k, l)$ is a so-called `Wythoff pair' or a `dual Wythoff pair', that is of the form $(\lfloor\phi n\rfloor, \lfloor\phi^2 n\rfloor)$ or 
$(\lceil\phi n\rceil, \lceil\phi^2n\rceil )$, $n\in \N$.
Thus, in Maharaja Nim we take the first Wythoff pair 
$(1, 2) = (\lfloor\phi\rfloor, \lfloor\phi^2\rfloor)$, whereas in the next 
section we study $(2,3)$-Maharaja Nim, that is we let $(k,l)$ take 
the values of 
the first dual Wythoff pair $(2,3) = (\lceil\phi \rceil, \lceil\phi^2\rceil )$.

\begin{figure}[ht]
\centering
\includegraphics[width=0.5\textwidth]{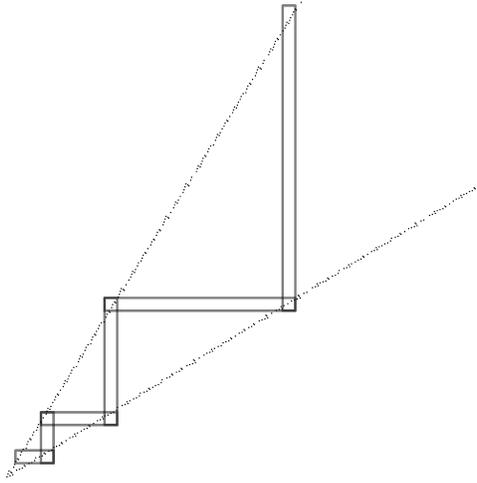}
\caption{A `telescope' with `focus' $O(1)$ and `reflectors' along 
the lines $\phi n$ and $n/\phi$ attempts to determine the outcome (P or N) 
of some position, $(x, y)$ at the top of the picture. As we demonstrate 
in Section \ref{Section:5} the method is successful for 
$(2,3)$-Maharaja Nim. (It gives the correct value 
for all extensions of Wythoff Nim with a finite non-terminating converging 
dictionary). The focus is kept sufficiently wide (a constant)
to provide correct translations in each step. 
The number of steps is linear in $\log (xy)$.}
\label{figure:12}
\end{figure}

\begin{Conj}\label{conjecture:1}
Let $k,l\in\N,$ $k<l$. Then each upper P-position $(x, y)$ 
of $(k, l)$-M satisfies $y = \phi x + O(1)$.
\end{Conj}
Does this conjecture hold for any game of the form `a 
finite number of moves adjoined to Wythoff Nim'?

Suppose that a given game $(k, l)$M has a finite (non-terminating) 
dictionary (as for Maharaja Nim) thus, hypothetically, 
providing an affirmative answer to Conjecture \ref{conjecture:1}. 
Suppose further that the dictionary \emph{converges}, 
that is, given an arbitrary string-position, we can, within the distance 
of a bounded number of bits, precisely determine when a new word starts. 
For this particular game, let us sketch a polynomial time algorithm which 
determines whether a given position $(x, y)$ (with $\frac{y}{x}$ approximately 
$\phi$) is P, see also Figure \ref{figure:12}. 
Suppose that we have computed an initial (sufficiently large) 
sequence of the bit-string.
We sketch the steps of the decision problem of $(k,l)$M are as follows:
\begin{itemize}
\item Back track $(x, y)$ via orthogonal reflections along the 
lines $\phi n$ and $n/\phi$. Here we do not need to use our dictionary, 
only to put marks at the precise locations of our reflecting points on 
the lines $\phi n$ and $n/\phi$. That is, we get a finite sequence of 
pairs of the form 
$$(x, \phi x), (x, x/\phi), (x/\phi^2, x/\phi),\ldots , (x/\phi^p, x/\phi^{p-1}),$$ 
some $p\in \N$.
\item When we have back tracked as far as to our initial bit-string, 
the `forward' translations can begin. Suppose that we know that the 
dictionary converges within $q$ (which is supposed to be 
much less than $x$ and $y$) bits and that the maximal length of 
a translate is $c\le q$ bits. 
\item Then it suffices to translate 
$< \phi q$ bits in each step. If the first left hand side word begins with, 
say the bit 
$\lfloor x/\phi^p\rfloor - \phi q \le b_1 \le \lfloor x/\phi^p\rfloor - q$ 
we may translate it and be assured to find another left hand side word 
beginning at a bit 
$\lfloor x/\phi^{p-1}\rfloor -\phi q \le b_2 \le \lfloor x/\phi^{p-1}\rfloor -q$ 
and so on. For the final computation of the value of $(x,y)$ it suffices to, 
given the left hand side word which contains $x$,
compute the P-positions in some area of size less than $ c\times c$ squares. 
(Alternatively, given a short dictionary, the list of P-positions 
corresponding to each word may be computed beforehand.) 
\item This procedure takes $p$ steps where $\phi^p$ is proportional to $x+y$.  
\end{itemize}

\section{The close relative $(2,3)$-Maharaja Nim has polynomial time complexity}\label{Section:5}
The game $(2, 3)$-Maharaja Nim, $(2, 3)$-M, is as Maharaja Nim except 
that, for this game, the Knight's jumps are of the form $(2,3)$ and $(3,2)$ 
(and not $(1,2)$ and $(2,1)$). In this section we let
$(a_1,b_1), (a_2,b_2), \ldots$ denote the upper P-positions of $(2,3)$-M, 
where $(a_i)$ is increasing. As we have remarked in Section \ref{Sec:4}  
analogs of Proposition \ref{rowcolumn} and (\ref{diag}) hold for $(2,3)$-M. 
Hence $(a_i)$ and $(b_i)$ are complementary.

\begin{figure}[ht]
\centering
\includegraphics[width=0.35\textwidth]{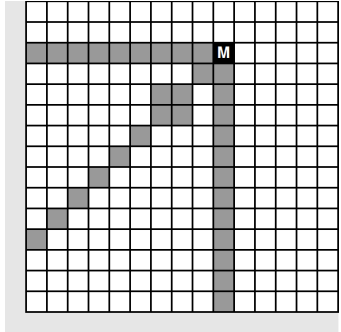}
\caption{The move options from a given position of $(2,3)$-Maharaja Nim.}\label{figure:13}
\end{figure}

\begin{figure}[ht!]
\centering
\includegraphics[width=0.55\textwidth]{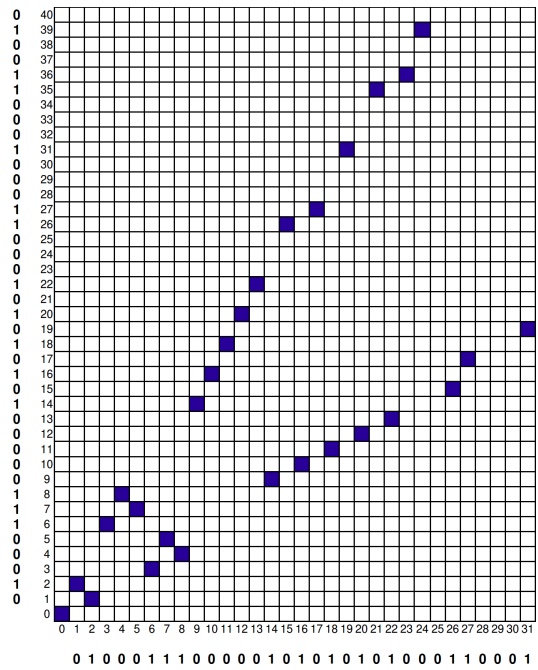}
\caption{The initial P-positions of $(2,3)$-Maharaja Nim together 
with its initial bit-string.}\label{figure:16}
\end{figure}

Since Lemma \ref{LemmaPerfect} does not hold for $(2,3)$-M, for the 
analysis of this game we use a relaxation of the approach in 
Section \ref{Sec:2}. As we saw at the end of that section, the crucial 
property for approximate linearity to hold is that the dictionary promised 
a sufficiently frequent reappearance of property (\ref{nth}). 
Hence, for a new left hand side word to be translated it is not necessary 
that we require a perfect sector as defined for $(1,2)$-Maharaja to be detected. 
It turns out that the condition (\ref{12n}) in Lemma \ref{LemmaPerfect} suffices for our purposes.
This almost corresponds to a perfect sector, by which we mean that at most a finite number 
of positions are deleted from a perfect sector. 
That is, the requirement is still that (\ref{nth}) and (\ref{12n}) hold simultaneously. 

Suppose that the initial P-positions up to column $a_n$ has been coded in a 
unique $(2,3)$-M bit-string, where 
as before, a `1' (`0') in th $i^{th}$ position denotes a lower (upper) 
P-position in column $i$. That is the read head is about to read the 
${a_n}^{th}$ bit in the string. As in Section \ref{Sec:2}, by 
symmetry of P-positions, a finite number of bits follow to 
the right of the read head's current position.
Then a (new) left hand side word $\omega\ne 1$ (the word `1' is translated to `0') 
is included to the dictionary if and only if the following two criteria 
are satisfied. Each one of the numbers $0, 1, \ldots , n-1$ is represented as 
the difference $b_i-a_i$ of the coordinates of precisely 
one of the first $n - 1$ upper P-positions \emph{and}  $b_n-a_n=n$. 

As usual, the translation of $\omega$ is computed and concatenated at the end 
of the bit-string. The next left hand side word begins by the ${a_n}^{th}$ column.

\subsection{$(2,3)$-Maharaja Nim's dictionary process}

Given a finite binary dictionary, we define unique \emph{non-prefix free} translations by the following rule. Suppose that the read head has finished one translation in the (infinite) binary string $x$ and starts reading at position $n$. Suppose further that it detects the left hand side entries $\omega_1,\ldots , \omega_k$ of the dictionary, reading from position $n$ and onwards, where $\omega_i$ is a prefix of  $\omega_{i+1}$ for all $i<k$, so that $\omega_k$ is not the prefix of any other entry. Then it accepts the translation of $\omega_k$ and it is unique if it exists.

Let us illustrate this definition by defining the following (very short) non-prefix free dictionary of $(2,3)$-M:

\begin{align}
0 &\rightarrow 10\label{0}\\
1 &\rightarrow 0\label{1}\\
01000 &\rightarrow 100011100\label{01000}\\
01010 &\rightarrow 10001100\label{01010}.
\end{align}

Since the bit `0' is a prefix of the words `01000' and `01010' we need 
some external rule to decide which translation to use in the construction 
of the bit-string.  Suppose that the next bit detected by the read head 
is `0'. Then the translation is as in (\ref{0}), except if the next 
four bits are either `1000' or `1010'. For these cases the 
translations are as in (\ref{01000}) and (\ref{01010}) respectively. 
Notice that, by these translation rules, by (\ref{0}) and (\ref{1}), \emph{any} bit-string 
has a (longer) translation, and therefore the construction cannot terminate. Before proving that 
this dictionary is correct, let us provide some initial translations.

Column-wise, the first non-terminal 
P-positions of $(2,3)$M are $(1, 2),$ $(2, 1), (3, 6), (4, 8), (5, 7), 
(6, 3), (7, 5)$ and $(8, 4)$. These P-positions 
correspond to the bit-string `01000111' on the $x$-axis and 
`0100011100000' on the $y$-axis, see Figure \ref{figure:16}. That is, we can assume that the first word to be translated 
starts in position $(9, 14)$, which corresponds to the first diagonal (of the form (\ref{diag})) in the $9^{th}$ column. It is clear the the initial 
interference between rows and columns has ended here so, to begin with 
the read heads position is `01000111\underline{0}0000'. Thus, the first 
four translations are of the type (\ref{0}) which produce the bit-string 
 `010001110000\underline{0}10101010'. Then a type (\ref{01010}) translation 
follows which produces `01000111000001010\underline{1}01010001100', and so on.

It is easy to see that, given a perfect 
sector to the right of the column of the read head's position, each 
translation in $(2,3)$-M's dictionary is correct. 
However, since there is no a priori guarantee that a new translation 
starts at a perfect sector, we need to exclude certain combinations of translations, thus preventing any $(2,3)$-type move to short-circuit two P-positions. The translations (\ref{0}) and (\ref{01000}) could potentially interfere with a succeeding translation but (\ref{1}) and (\ref{01010}) cannot. Precisely, if the word `0' were followed by a `1' and then any of the words beginning with `0', or if the word `01000' were followed by a `0', then the translation rules would be wrong, because of a $(2,3)$-type ``short-circuit" of P-positions. These are all cases that we need to exclude. Let us begin to rule out the latter case.\\

\noindent {Claim 1:} If the left hand side word 
`01000' is detected by the read head, then it is succeeded by the left hand side word `1'.\\

Suppose, on the contrary, that the read 
head reads the left hand side word `01000' followed by a `0'. 
This string, `010000', which we say is part of our \emph{original} string, must have been translated from the left hand side words  
`$x$',`0','1','1','1' (in this order and where $x$ is the left hand side word in either 
(\ref{0}), (\ref{01000}) or (\ref{01010})). But the string `00111' 
only appears as a translation in (\ref{01000}). Further, the string `01000111' 
is forced since `11000111' cannot appear, but it cannot be that the read 
head detected the first five bits `01000' 
as the word in (\ref{01000}), since it would have translated to `100011100000' which does not include the original string '01000' in the right place. 
Thus, to prevent this, preceding the pattern `01000111' there must have been either  `0100', `01', or `0101'. It follows that either of the strings 
\begin{align}
& 010001000111,\label{first}\\ 
& 0101000111, \text{ or}\label{second}\\ 
& 010101000111\label{third} 
\end{align}
must have been read in the stage before the original string, \emph{and} where a new left hand side word starts from the first `0'. 
But then, for the case (\ref{first}), 
this translates to the original pattern `1000111000101010000', which 
forces that a left hand side word starts after the consecutive words `1',`1',`1', that is `01010' will be detected as a word and so the word `01000' would not have been read in the original string, which contradicts our initial assumption. 

For the latter cases (\ref{second}) and (\ref{third}), we get the translates `100011001010000' and `100011000101010000' respectively, both which may be treated in analogy to the first case, but here it is forced that new left hand side words start after the consecutive words `1',`1' respectively. 
\hfill $\Box$\\ 

\noindent {Claim 2:} Any sequence of left hand side words beginning with `0', `1' and then some pattern beginning with a `0' is impossible.\\

We are here concerned with that the read head  detects any sequence of left 
hand side words beginning with `0', `1' and then some sequence 
`$0xy$', where $x$ and $y$ represent two bits. 
It is immediate by the translation rules that we may exclude the cases where 
$xy$ represents `00' or `10'. Namely, for these two cases, by 
the `prefix-rule' of choosing the longest left hand side word in 
the dictionary, we would rather have used one of the translations 
in (\ref{01000}) or (\ref{01010}). Also, the case where $xy$ is `11' may be excluded since 
the string `01011' does not appear in any combination of the right hand side translates. Thus, it only remains to analyze the case where the two bits are `01'. 
That is, we want to exclude the pattern `01001'. By looking at the translations it is obvious that the string `1001' must have been translated from the left 
hand side words `0', `1' and then a word beginning with a `0'. This means that 
precisely the pattern which we want to exclude has appeared in a 
previous translation (and thereby also short-circuiting two P-positions in columns strictly to the left of the current position). Thus (using Figure \ref{figure:16} as a base case) strong induction resolves this case.

\subsection{Polynomiality}
We have proved that the dictionary in (\ref{0}) to (\ref{01010}) is correct 
and thereby also that the P-positions of 
$(2,3)$-Maharaja Nim lie within a bounded distance of either the 
`line' $\phi n$ or $\phi^{-1}n$. Next, we will demonstrate that 
this dictionary gives a polynomial strategy, as outlined 
in Section \ref{Sec:4}. 
For this, it suffices to prove that, given an arbitrary position in 
the infinite bit-string, by a search within a bounded number of bits 
we can determine which one of the four given translations is correct. 

If the read head reads the pattern `11' then, by the left hand side words 
in the dictionary and in particular (\ref{1}), 
we can conclude that a new word starts by the first `1'. 
Hence we assume that no two consecutive `1's are detected. 
By analyzing the translations in the dictionary one can see that at most five consecutive `0's can appear. Therefore, we may assume that the read 
head reads the pattern `010' within a bounded distance, 
which by previous arguments mean either `01000' or `01010'. Both these strings are detected as words, unless the preceding pattern ends 
with `0100', `01' or `0101'. Hence one needs to investigate the following 
six ambiguous strings:

\begin{enumerate}[(a)]
\item 010001000,
\item 0101000,
\item 010101000,
\item 010001010,
\item 0101010,
\item 010101010.
\end{enumerate}
The pattern `10001000' in (a) cannot have been translated from the string `011011'. This follows by viewing the possible combinations of right hand side translates. Hence, the combination of translations comes from first `0',  `1' and  `1' and then `01000' or `01010'. But these combinations are also impossible since they both enforce the impossible pattern ``1101''. Hence (a) cannot appear.

The string in (b) must have been translated from `0,', `0', `1', `1' 
which, by (\ref{01000}) and (\ref{01010}) and since all 
translates end with a `0', implies that the three preceding bits must have been `010'. Hence, we can extend the pattern to be translated to `0100011'. It 
is given that the prefix `01000' of this string cannot be 
detected as a left hand side word. Therefore, the translation 
of `0100011' must be `10010101000' which has the prefix
\begin{align}\label{1001}
\text{`1001'}. 
\end{align}
But, by the left hand side words in the dictionary, 
any string containing (\ref{1001}) must converge between the two `0's. 
Hence a new word must start as `01010' followed by `1', `0', `0',$\ldots$. 
Notice that (c) has this string as a suffix and hence 
it may also be included in the argument. Also, by (\ref{1001}) and by the argument in (a), in any attempt to disprove convergence (d) must be preceded by the pattern `01', but then again, we may analyze (d) as (b).

We are left with the strings (e) and (f). 
Since a repetition of more than five consecutive patterns `01' implies that 
more than five consecutive 0s has been translated, which is impossible, we 
may assume that the repetitions of `01' in (f) has been preceded by 
either of the patterns `10' or `00' (`11' is already ruled out). Again, 
the first case leads to (\ref{1001}). Notice that (e) can also 
be included in this argument. For the second case, notice that any 
string beginning with `00001' converges after the three first `0's, that is 
a new word must begin with `01', so it suffices to study the string 
`1000101010', which (since the pattern `11' is excluded) has been 
treated already in (d).
 
We have proved that, given an arbitrary position in the bit-string, at most 
a bounded number of preceding bits need to be searched in order to find the 
correct translation. By Section \ref{Sec:4} this convergence gives a polynomial 
time winning strategy of $(2,3)$-Maharaja Nim.
\appendix
\numberwithin{equation}{section}
\section{Code}
\subsection{The Maple code corresponding to Figure \ref{figure:2}}\label{A:1}
The below code includes the P-positions of both Wythoff Nim and Maharaja Nim in one and the the same diagram. 
\small{
\begin{verbatim}
restart: with(plots): with(plottools):

N:=50;

theLine1:=CURVES([[0.0,0.0], [evalf(N), evalf(N*(1+sqrt(5))/2)]]): 
theLine2:=CURVES([[0.0,0.0], [evalf(N*(1+sqrt(5))/2), evalf(N)]]):

#Compute the P-positions of Wythoff Nim and store as a list of squares.
#0=Not yet computed, 1=P, 2=N.
for i from 0 to N do for j from 0 to N do A[i,j]:=0: od: od:
for i from 0 to N do for j from 0 to N do if A[i,j]=0 then A[i,j]:=1: 
for k to N do A[i+k, j]:=2: A[i+k,j+k]:=2: A[i,j+k]:=2: od: fi: od: od:
rectListW:=[]: for i from 0 to N do for j from 0 to N do if A[i,j]=1 
then rectListW:=[op(rectListW), [[i,j],[i,j+1],[i+1,j+1],[i+1,j]]]: fi: 
od: od:

#Draw the P-positions and the two lines with slopes the golden ratio:
display(polygonplot(rectListW, color=red), theLine1, theLine2, axes=none, 
scaling=constrained, view=[0..N, 0..N]);

#Compute the P-positions of Maharaja Nim:
for i from 0 to N do for j from 0 to N do A[i,j]:=0: od: od:
for i from 0 to N do for j from 0 to N do if A[i,j]=0 then A[i,j]:=1:
A[i+1,j+2]:=2:
A[i+2,j+1]:=2:
for k to N do A[i+k, j]:=2: A[i+k,j+k]:=2: A[i,j+k]:=2: od: fi: od: od:
rectListM:=[]: for i from 0 to N do 
for j from 0 to N do if A[i,j]=1 then rectListM:=[op(rectListM), 
[[i+0.2,j+0.2],[i+0.2,j+0.8],[i+0.8,j+.8],[i+0.8,j+0.2]]]: fi: od: od:

display(polygonplot(rectListM, color=blue), axes=none, scaling=constrained);
display(polygonplot(rectListM, color=blue), 
polygonplot(rectListW, color=red), theLine1, theLine2, axes=none, 
scaling=constrained, view=[0..N, 0..N]);
\end{verbatim}}

\subsection{The Maple code corresponding to Maharaja Nim's dictionary.}\label{A.2} 
The following code explores whether the first 9 words in Maharaja Nim's dictionary suffices.
{\small
\begin{verbatim}

dictionary:={[1], [0,1], [0,0,1,0,0], [0,0,1,0,1,1,0], [0,0,1,1,0],
[0,0,0,1,0,0],[0,0,0,0,1,0,0,1,0], [0,0,0,0,0,1,0,0], [0,0,1,1,1,0]};

translation:=table([[1]=[0], [0,1]=[1,0,0], [0,0,1,0,0]=[1,0,0,1,0,1,1,0,0],
[0,0,1,0,1,1,0]=[1,0,0,1,0,0,1,1,0,0,0],[0,0,1,1,0]=[1,0,0,1,0,1,0,0], 
[0,0,0,1,0,0]=[1,0,0,1,0,1,1,0,1,0,0], 
[0,0,0,0,1,0,0,1,0]=[1,0,0,1,0,0,1,1,1,1,0,0,0,1,0,0], 
[0,0,0,0,0,1,0,0]=[1,0,0,1,0,1,1,0,0,1,1,1,0,0,0], 
[0,0,1,1,1,0]=[1,0,0,1,0,0,1,0,0]]);

theString:=[0,0,1,0,0]: reader:=0:
for times to 12000 do foundWord:=false: 
for i to 9 do if not foundWord then theWord:=theString[reader+1..reader+i]: 
if member(theWord, dictionary) then foundWord:=true: 
theString:=[op(theString), op(translation[theWord])]: 
reader:=reader+i: fi: fi: od: if not foundWord 
then print(reader, theString[reader+1..reader+20]): fi: 
if times mod 100 = 0 then print(times, nops(theString)): fi: od:

\end{verbatim}}
\section{An undecidable dictionary process}

First we describe a known undecidable problem.

Suppose we have an alphabet consisting of one special symbol S which acts as "space" or "stop" symbol, and a finite number of other symbols denoted A, B, C,$\ldots$

We start from a "multiplication table" that describes an operation $x*y = z$, where $x$ and $y$ are arbitrary symbols from the alphabet, and $z$ is a symbol other than S.

For instance, the table may look like
\begin{center}
  \begin{tabular}{ l | c c c c }
    & S & A & B & C \\\hline
S & A & B & B & B \\
A & C & B & A & A \\ 
B & C & C & C & C \\
C & A & A & B & A \\   
  \end{tabular}
\end{center}

Given such a table, we form a triangular pattern of symbols consisting of rows starting and ending with $S$, and where the other symbols are obtained by "multiplying" the two symbols above it.

The table in the example gives

\begin{center}
  \begin{tabular}{c}
   S\\
   S S\\
   S A S\\
   S B C S\\
   S B C A S\\
   S B C A C S\\
   S B C A A A S\\
   S B C A B B C S\\
  \end{tabular}
\end{center}
and so on.

Naively, one would like to understand how this pattern behaves by looking at the table. The hope of general understanding of this kind is shattered by the fact that the behavior of the pattern can simulate any given Turing machine. It follows that a number of simple questions are generally algorithmically undecidable. We mention a few such questions which are easily seen to be "equivalent".

If we are given two multiplication tables, do they produce the same pattern of symbols or not? We can examine the tables and find the entries where they differ. If only we can decide whether any of these entries is ever going to be used, we are done. If we fill in the entries where the tables differ with a "new" symbol Z, then in turn the problem becomes equivalent to deciding whether or not a certain symbol of the alphabet is ever going to occur in the pattern. This question in turn is equivalent to deciding whether a partial multiplication table (one with empty places) is "consistent" in the sense of determining a pattern.

In the example above, it is straightforward to see that the table entries B$*$S and S$*$C are never going to occur, but that all other entries do. Thus if we change the two entries B$*$S and S$*$C to something else, the pattern will still be the same, while if we change any other entry, the resulting pattern will be different. But in general, answering such questions may be as difficult as any mathematical problem. For instance, it is possible to "program" a multiplication table to look for counterexamples to the Goldbach conjecture, so that a certain symbol of the alphabet occurs if and only if there is a large even number which is not the sum of two primes.

Now consider a different type of process. Here we may without loss of generality assume the alphabet to be $\{0,1\}$. Suppose we have a given starting string A, and a "dictionary" consisting of "translations" of the form $x \rightarrow y$, where $x$ and $y$ are binary strings. The dictionary is a finite set of such translations $$x_1 \rightarrow y_1, \ldots , x_n \rightarrow y_n,$$ and to avoid ambiguity, it is required that no $x_i$ is a prefix of any other. (Although generalized dictionaries such as those in Section \ref{Section:5} are also undecidable for the same reasons as explained here).

A "reader" starts at the left endpoint of the string A, and "reads" until it finds a word $x_i$. Then a "writer", initially at the right endpoint of A, writes the translation $y_i$ and concatenates it to the right of A. The reader then continues from where it was interrupted, and reads until it finds the next word etc.

The process may either go on forever, or get stuck by the reader reaching the right endpoint of the string without finding a word in the dictionary. The analogous questions may be asked about this process. Does it terminate or not? Is a certain word ever read? We will show that the "multiplication process" can be encoded as a "dictionary process", thereby showing that in general, the fundamental questions about the dictionary process are undecidable.

Suppose therefore that we are given a multiplication table. We will construct a dictionary that mimics the pattern of symbols arising from the given multiplication table.

First we introduce "metasymbols" that are binary strings representing the symbols of the alphabet. The starting string is going to be SS, and there is one dictionary entry for each entry of the multiplication table. If the table contains, for instance, A$*$B $=$ C, then there is a translation rule AB $\rightarrow$ CC.

The idea is that instead of writing a C, the writer writes CC. Eventually, the reader will read the second "half" of the previous symbol together with the first C, and then the second C together with the first "half" of the next symbol, while the writer produces the corresponding products. To achieve this, the symbol S needs special treatment. We therefore also use the translation rules SS $\rightarrow $S(S$*$S)(S$*$S)S, and for each other symbol A, SA $\rightarrow$ S(S$*$A)(S$*$A) and AS $\rightarrow$ (A$*$S)(A$*$S)S.

In the example above, we get the dictionary\\

\noindent SS $\rightarrow$ SAAS\\
SA $\rightarrow$ SBB\\
SB $\rightarrow$ SBB\\
SC $\rightarrow$ SBB\\
AS $\rightarrow$ CCS\\
AA $\rightarrow$ BB\\
AB $\rightarrow$ AA\\
AC $\rightarrow$ AA\\
BS $\rightarrow$ CCS\\
BA $\rightarrow$ CC\\
BB $\rightarrow$ CC\\
BC $\rightarrow$ CC\\
CS $\rightarrow$ AAS\\
CA $\rightarrow$ AA\\
CB $\rightarrow$ BB\\
CC $\rightarrow$ AA\\

Starting from the string SS, this produces (and here we have introduced some spacing just to increase readability)\\

SS SAAS SBBCCS SBBCCAAS SBBCCAACCS SBBCCAAAAAAS 
$\ldots$\\

By leaving out certain rows of the dictionary, we may mimic a partial multiplication table. Therefore the question whether the dictionary process terminates is algorithmically undecidable.

\end{document}